\documentclass[12pt]{amsart}

\usepackage[dvips]{graphicx}

\usepackage[centertags]{amsmath}

\usepackage{amsfonts,amssymb,latexsym}

\usepackage{fullpage}

\DeclareMathAlphabet{\mathbcal}{U}{eur}{b}{n}

\newcommand{\Dt}{\Delta}

\newcommand{\bbE}{{\mathbb E}}

\newcommand{\bpf}[1][Proof]{{ {\sc #1: }}}

\newcommand{\cN}{{\mathcal N}}

\newcommand{\epf}{{$\quad\Box$}}

\newcommand{\real}{{\RR}}
\newcommand{\bx}{\overline{X}}
\newcommand{\kk}{\kappa}
\renewcommand{\Sigma}{g}

\newcommand{\ks}{k^{\star}}
\newcommand{\by}{{\overline {y}}}
\newcommand{\ddelta}{\varepsilon}

\newcommand{\bbR}{\mathbb{R}}
\newcommand{\RR}{\bbR}

\newcommand{\la}{\langle}
\newcommand{\ra}{\rangle}

\newcommand{\bX}{{\overline{X}}}

\newcommand{\FFrac}[2]{{\textstyle \frac{#1}{#2} }}

\newcommand{\E}{\mathbb{E}}
\newcommand{\Prob}{\mathbb{P}}
\newcommand{\ZZ}{\mathbb{Z}}

\newcommand{\calO}{{\mathcal O}}

\newtheorem{theorem}{Theorem}
\newtheorem{lemma}{Lemma}[section]

\newtheorem{corollary}[lemma]{Corollary}

\newtheorem{assum}[lemma]{Assumption}

\newcommand{\gm}{\gamma^{-}}

\newcommand{\bg}{\bar{\gamma}}
\newcommand{\xit}{\tilde{\xi}}
\newcommand{\Mt}{\tilde{M}}

\newcommand{\ONE}{\mathbf{1}}

\newcommand{\EE}{\mathbb{E}}
\newcommand{\NN}{\mathbb{N}}
\newcommand{\PP}{\mathbb{P}}

\newcommand{\eqdef}{\:{\overset{\mbox{\tiny def}}{=}}\:}

\newcommand{\cB}{\mathcal{B}}

\newcommand{\cF}{\mathcal{F}}

\newcommand{\cH}{{\mathcal H}}
\newcommand{\cG}{\mathcal{G}}

\newcommand{\ta}{{\tilde {\alpha}}}
\newcommand{\tb}{{\tilde {\beta}}}
\newcommand{\dm}{{\delta^{-}}}
\newcommand{\dpp}{{\delta^{+}}}
\newcommand{\xs}{x^{\star}}

\newcommand{\ps}{p^{\star}}
\newcommand{\qs}{q^{\star}}

\begin{document}

\title[An Adaptive Euler-Maruyama Scheme For SDEs: Convergence
and Stability]{An Adaptive Euler-Maruyama Scheme For SDEs: Convergence and
Stability}
\author[H. Lamba, J.C. Mattingly \and  A.M. Stuart]{H. Lamba$^{1}$,
J.C. Mattingly$^{2}$ \and  A.M. Stuart$^{3}$}

\footnotetext[1]{Department of Mathematical Sciences, George Mason University,
Fairfax, VA22030, USA. $${\tt hlamba@gmu.edu}$$}
\footnotetext[2]{Department of Mathematics, Duke University,
  Durham, NC 27708, USA. $${\tt jonm@math.duke.edu}$$}
 \footnotetext[3]{Mathematics Institute, Warwick University, Coventry, CV4 7AL, England. $${\tt andrew.m.stuart@warwick.ac.uk}$$}

\begin{abstract}

The understanding of adaptive algorithms for SDEs is an open
area where many issues related to both convergence and stability (long time
behaviour) of algorithms are unresolved. This paper considers a very simple
adaptive algorithm, based on controlling only the drift component
of a time-step. Both convergence and stability are studied.

The primary issue in the convergence analysis is that the
adaptive method does not necessarily drive the time-steps to zero
with the user-input tolerance. This possibility must be quantified and
shown to have low probability.
  
The primary issue in the stability analysis is ergodicity.
  It is assumed that the noise is
  non-degenerate, so that the diffusion process is elliptic, and the
  drift is assumed to satisfy a coercivity condition. The SDE is
  then geometrically ergodic (averages converge to statistical equilibrium
  exponentially quickly). If the drift is not linearly
  bounded then explicit fixed time-step approximations, such as the
  Euler-Maruyama scheme, may fail to be ergodic.  In this work, it is
  shown that the simple adaptive time-stepping strategy cures this
  problem.  In addition to proving ergodicity,
  an exponential moment bound is also proved,
generalizing a result known to hold for the SDE itself.

\vspace{0.2in}

\end{abstract}
\maketitle

\section{Introduction}
\setcounter{equation}{0}

In this paper, we study the numerical solution of the It\^{o} SDE
\begin{equation}
   dx(t) = f(x(t)) dt + g(x(t)) dW(t), \quad x(0) = X
 \label{eq:sde}
\end{equation}
by means of an adaptive time-stepping algorithm.
Here $x(t) \in \real^m$ for each $t$ and $W(t)$ is a $d-$dimensional
Brownian motion. Thus $f: \real^m \to \real^m$ and $g:\real^m \to
\real^{m \times d}.$ For simplicity we assume that the initial condition is
deterministic. Throughout
$|\cdot|$ is used to denote either the Euclidean vector norm or
the Frobenius (trace) matrix norm as appropriate.
We assume throughout that $f,g$ are $C^{2}$.
Further structural assumptions will be made where needed. 
The basic adaptive mechanism we study is detailed at the start of the
next section.  It is a simple adaptive algorithm, prototypical of a whole class
of methods for the adaptive integration of SDEs.  Our aim is twofold. First we 
show convergence, as the user-input tolerance $\tau$ tends to zero;
this is a non-trivial
exercise because the adaptive strategy does not imply that the
time-steps taken tend to zero with the tolerance everywhere in phase space.
Secondly we show that the methods
have a variety of desirable properties for the long-time integration of
ergodic SDEs, including preservation of ergodicity and exponential
moment bounds.

The adaptive method controls the time-step of a forward Euler drift step
so that it deviates only slightly from a backward Euler step. 
This not only controls an estimate of the contribution to the
time-stepping error from the drift step, but also allows the analysis of 
stability (large time) properties for implicit backward Euler
methods to be employed in the explicit adaptive method. 
Numerical experiments suggest that both the convergence and stability
analyses extend to a number of more sophisticated methods which control
different error measures; some of these experiments are
reported below.

It is of interest to discuss our work in the context of
a sequence of interesting papers which study the optimality
of adaptive schemes for SDEs, using various different error
measures \cite{ritter, ritter2, ritter3, MG2}. For many of these error measures,
which are quite natural in practice,
the asymptotically optimal adaptive schemes are based solely on the diffusion.
This is essentially because it is the diffusion term which dominates
the (lack of) regularity in paths and this regularity in  turn
dominates error measures.
Why then, have we concentrated on methods which adapt only
on the drift?  The reason for this is that, as mentioned above,
such methods are advantageous for long-time integration.
In practice, we anticipate that
error controls based on both drift and diffusion could combine
the advantages of the asymptotically optimal schemes with the
enhanced stability/ergodicity of schemes which control based
on the drift.

In order to prove a strong mean-square
convergence result for this algorithm, it is first necessary to obtain a
suitable upper bound on the sequence of timesteps used. These bounds
mimic those used in the convergence proofs for adaptive ODE solvers
\cite{S, HL, L} and require that the numerical solution does not enter
neighbourhoods of points where the local error estimate vanishes. (Requiring
that these neighbourhoods are small excludes some simple drift
vector fields, such as constants. In practice, we would anticipate controlling
on both the drift and the diffusion, minimizing this issue). An
essential part of the analysis is a proof that the contribution to the
mean-square error from paths that violate this condition is suitably
small.

Adaptivity is widely used in the solution of ordinary differential
equations (ODEs) in an attempt to optimize effort expended per 
unit of accuracy. The adaptation strategy can be viewed
heuristically as a fixed time-step algorithm applied to
a time re-scaled differential equation \cite{Gr} and it
is of interest to study convergence of the algorithms
as the tolerance employed to control adaptation is
reduced to zero \cite{HL}. However adaptation also
confers stability on algorithms constructed from
explicit time-integrators, resulting in better qualitative behavior
than for fixed time-step counter-parts. This viewpoint
was articulated explicitly in \cite{Sanz} and subsequently pursued
in \cite{AGH}, \cite{h&s} and \cite{SH2} for example. In
particular the reference \cite{SH2} studies the effect of
time-discretization on dissipative structures such as those
highlighted in \cite{Ha,Te}. It is shown that certain adaptive
strategies have the desirable property of constraining time-steps
of explicit integrators so that the resulting solution update
differs in a controlled way from an implicit method. Since
many implicit methods have desirable stability properties
(see \cite{DV} and \cite{SH}, Chapter 5) this viewpoint can
be used to facilitate analysis of the stability of adaptive
algorithms \cite{SH2}.

In \cite{msh}, stochastic differential equations (SDEs)
with additive noise and vector fields
satisfying the dissipativity structures of \cite{Ha,Te} are
studied. There, and in \cite{RT, Talnew, Talnew2}, it is shown that explicit
time-integrators such as Euler-Maruyama may fail to be ergodic even when the
underlying SDE is geometrically ergodic. The 
reason is that the (mean) dissipativity induced by the
drift is lost under time-discretization. Since this is
exactly the issue arising for explicit integration of dissipative ODEs, 
and since this issue can be resolved in that context by means of
adaptation, it is natural to study how such adaptive methods
impact the ergodicity of explicit methods for SDEs.
In recent years, the numerical solution of SDEs with gradient
drift vector fields has been used as the proposal for an MCMC method for sampling
from a prescribed density, known only up to a multiplicative
constant -- a technique referred to as Metropolis-adjusted Langevin 
algorithm \cite{CR}.
In this context, it is very desirable that the time discretization
inherit ergodicity. The adaptive scheme proposed here is an
approach to ensuring this.
In this sense our work complements a number of recent papers concerned
with constructing approximation schemes which are ergodic in situations where
the standard fixed step Euler-Maruyama scheme fails to be: in
\cite{RT} a Metropolis-Hastings rejection criterion is used to
enforce ergodicity; in \cite{hansen, s&t} local linearization is used;
in \cite{msh} implicit methods are used. 
Although the adaptive method that we analyze here
is proved to be convergent on finite
time intervals, it would also be of interest
to extend the work of Talay \cite{Tal}, concerned with convergence proofs
for invariant measure under time-discretization, to the adaptive
time-step setting considered here.

In Section \ref{sec:Algorithm}, we introduce the adaptive
algorithm, together with some notation. 
In Section~\ref{sec:res} the finite time convergence result for the
adaptive method is stated. The proof is given in Section~\ref{sec:alg}
and proceeds by extending the
fixed step proof given in \cite{HMS}; 
the extension is non-trivial because the adaptivity does not force
time-steps to zero
with the tolerance in all parts of the phase space.
In Section \ref{sec:mainResults}, we state the main results of the paper on the
stability of the adaptive method. All results are proved under
the dissipativity condition
\begin{equation}
\exists \alpha, \beta \in (0,\infty):
 \la f(x),x \ra \le \alpha-\beta |x|^2 \quad \forall x \in \real^m,
\label{eq:coerce} 
\end{equation}
where $\la \cdot, \cdot \ra$ is the inner-product inducing the Euclidean norm,
as well as a boundedness and invertibility condition on the diffusion matrix
$g$.  The results proven include ergodicity
and an exponential moment bound; all mimic known results about the
SDE itself under \eqref{eq:coerce}. 
Section \ref{sec:APriori} starts with a number of {\em a priori} 
estimates for the adaptive scheme of
Section \ref{sec:Algorithm}, and proceeds to proofs of the stability 
stated in Section \ref{sec:mainResults}.  
Numerical results studying both convergence and ergodicity are
presented in Sections \ref{sec:numerics}--\ref{sec:9}. Some concluding
remarks and generalizations are given in Section \ref{sec:conc}.

\section{Algorithm}\label{sec:Algorithm}
\setcounter{equation}{0}

The adaptive algorithm for \eqref{eq:sde} is as follows:
\begin{eqnarray}
\label{eq:methoda}
\begin{array}{ccc}
k_{n} &=& \cG(x_n,k_{n-1}), \quad k_{-1}=K\\
x_{n+1} &=& \cH(x_n,\Dt_{n})+\sqrt \Dt_{n} \Sigma(x_n)\eta_{n+1},
\quad x_0=X,\notag
\end{array}
\end{eqnarray}
where $\Dt_n=2^{-k_n}\Delta_{max}$. Here 
\begin{align*}
  \cH(x,t)&=x+tf(x)
\end{align*}
and
\begin{align*}
  \cG(x,l)&=\min\{k \in \ZZ^+:
  |f(\cH(x,2^{-k}\Delta_{max}))-f(x)| \le \tau \; \& \; k \ge l-1\}.
\end{align*}
The random variables $\eta_j \in \real^d$ form
an i.i.d. sequence distributed as $\cN(0,I).$
The parameter $K$ defines the initial time-step and $\tau>0$ the tolerance.
Note that the algorithm defines a Markov chain for $(x_n,k_{n-1})$ on 
$\RR^d \times \ZZ^+.$

We may write
\begin{eqnarray}
\label{eq:ratz}
\begin{array}{ccc}
\xs_n & = & x_n+\Delta_n f(x_n),\\
x_{n+1} &=& \xs_n+\sqrt \Dt_{n}\Sigma(x_n)\eta_{n+1}.
\end{array}
\end{eqnarray}
If $ K \in \ZZ^+$ then $k_n \in \ZZ^+$
and the error control enforces the condition
$$\Dt_n \le \min\{2\Dt_{n-1}, \Delta_{\max}\},$$
where $\Delta_{max}$ is the fixed maximum time-step.
Furthermore we have
$$|f(\xs_n)-f(x_n)| \le \tau.$$
In the absence of noise, this implies that
the difference between an Euler approximation at the
next time-step, and an explicit second order approximation, is of
size ${\calO}(\Delta_n\tau).$ In the presence of noise, it imposes
a similar restriction on the means.  
As mentioned in the introduction, in practice we would anticipate combining this
drift error control with others tuned to the diffusion.

\subsection{Notation}

The most important notation conceptually is concerned with making
relationships between the numerical approximations at discrete steps,
and the true solution  at certain points in time. To do this
we define ${\mathcal F}_n$ to be the sigma-algebra generated by $n$ steps of
the Markov chain for $(x_n,k_{n-1}).$ Let
\begin{equation*}
  t_n=t_{n-1}+\Dt_{n-1}, \quad t_0=0,  
\end{equation*}
$\delta>0$ and define the stopping times
$N_j$ by $N_0=0$ and, for $j \ge 1$,
\begin{equation*}
  N_j=\inf_{n \ge 0}\{n:t_n \ge  \delta+t_{N_{j-1}}\}.  
\end{equation*}
Where the dependence on $\delta$ is important we will
write $N_j(\delta).$
It is natural to examine the approximate process
at these stopping times since they are spaced approximately
at fixed times in the time variable $t.$ Theorem \ref{t:stf} in
Section \ref{sec:mainResults} shows that these stopping times are
almost surely finite, under the dissipativity condition \eqref{eq:coerce}.
Notice that
\begin{equation*}
  \delta^{-}:=\delta \le t_{N_j}-t_{N_{j-1}} \le \delta + 
  \Delta_{max}:=\delta^{+}.  
\end{equation*}
When considering strong convergence results
it is necessary to interpret $\sqrt \Delta _n \eta_{n+1}$ in the adaptive
algorithm as the Brownian increment $W(t_{n+1})-W(t_n).$

We let
\begin{align*}
y_j=x_{N_j+1} \quad\mbox{and}\quad l_j=k_{N_j}.  
\end{align*}
The Markov chain for $\{y_j, l_j\}$ will be important in our 
study of long time behaviour and we will prove that it is ergodic.
Let $\cG_{j}=\cF_{N_j}$, the filtration of events up to
the $j^{th}$ stopping time.

It is convenient to define two continuous time
interpolants of the numerical solution. We set
\begin{align}
X(t)&= x_n, \;\;t \in [t_n,t_{n+1}),\\
\bx(t)&= X+\int_{0}^t f(X(s))ds+\int_0^t \Sigma(X(s))dW(s).
\label{eq:bx1}
\end{align}
Hence, for $t \in [t_n,t_{n+1})$
\begin{align}
\bx(t) &= x_n+(t-t_n)f(x_n)+\Sigma(x_n)[W(t)-W(t_n)]\\
&=(1-\alpha_n (t))x_n+\alpha_n (t) \xs_n+\Sigma(x_n)[W(t)-W(t_n)]
\label{eq:bx2}
\end{align} 
for $\alpha_n (t)=(t-t_n)/(t_{n+1}-t_n) \in [0,1).$

It is sometimes important to know the smallest step-size
beneath which the error control is always satisfied, at a given point $x$.
Hence we define
\begin{align*}
  \ks(x)= \min\{k \in \ZZ^+:
  |f(\cH(x,2^{-l}\Delta_{max}))-f(x)| \le \tau \; \forall l \ge k\}
\quad\mbox{and}\quad \ks(B)=  \sup_{x \in B}
  \ks(x), 
\end{align*}
noting that, by continuity of $f$, $\ks(B)$ is finite if $B$ is bounded.

Because of the boundedness of $\Sigma$ we deduce that there are
functions $\sigma(x)$ and constants $\sigma, a>0$ such that, for $\eta$
distributed as $\eta_1$ and independent of $x$, 
\begin{equation*}
  \bbE |\Sigma(x)\eta|^2:= \sigma^2(x) \le \sigma^2,\quad\mbox{ and }\quad
\bbE \big||\Sigma(x)\eta|^2-\sigma^2(x)|^2 \big| \le  a\sigma^4.
\end{equation*}
The following definitions will be useful:
\begin{equation*}
\ta=\alpha+\frac12 \tau, \quad \tb=\beta-\frac12 \tau,
\quad \beta_n=\frac1{1+2\tb\Dt_{n}}, \quad \bg=1+\tb\Delta_{max}. 
\end{equation*}
We will always assume that $\tau$ is
chosen small enough so that $\tb>0$.  
The constants $\gamma^{-}$ is chosen so that 
$$(1+t)^{-1} \le  (1-\gm t) \le e^{-\gm t} \quad \forall t \in [0,2\tb \Delta_{max}].$$

\section{Convergence Result}\label{sec:res}

We start by discussing the error control mechanism.
We define $F_1(u)$ by
$$F_1(u) = df(u)f(u),$$ 
the function $F_2(u,h)$ by
$$F_2(u,h):=h^{-1}\Bigl(f(u+hf(u))-f(u)-hF_1(u)\Bigr)$$
and $E(u,h)$ by
$$E(u,h)=f(u+hf(u))-f(u).$$
Now, since $f \in C^2$, Taylor series expansion gives
\begin{equation} 
  E(u,h) =  h[F_1(u)+hF_2(u,h)]
\label{eq:form}
\end{equation}
where $F_1, F_2$ are defined above.
Note that the error control forces a time-step so that
the norm of $E(x_n,\Delta_n)$ is of order ${\calO}(\tau).$
Estimating the implications of this for the time-step $\Delta_n$ forms
the heart of the convergence proof below.

In order to state
the assumptions required for the convergence result we define, for
$R,\epsilon \geq 0,$ the sets
\begin{equation*}
\Psi(\epsilon) = \{u\in \RR^m : |F_1(u)| \leq \epsilon\},\;\;
 B_{R} = \{u \in \RR^m : |u| \leq R \}\;\;
{\rm and }\;  B_{R,\epsilon} = B_R \setminus \Psi(\epsilon)  
\end{equation*}
and introduce the constant $K_R = \sup_{u \in B_R,h\in [0, \Delta
  t_{\rm max}]}|F_2(u,h)|.$
Now define the following:
\begin{alignat*}{3}
\sigma_R &:= \inf\{ t \ge 0 : |\bX(t)| \ge R \}, &
  \rho_R &:= \inf\{ t \ge 0 : |x(t)| \ge R \},& \theta_R &:= \sigma_R  \wedge \rho_R \\
  \sigma_\epsilon &:= \inf\{ t \ge 0 : |F_1(\bX(t))| \leq \epsilon \},
  &
  \rho_\epsilon &:= \inf\{ t \ge 0 : |F_1(x(t))| \leq 2\epsilon
  \},&\theta_\epsilon &:= 
  \sigma_\epsilon  \wedge \rho_\epsilon \\
  \sigma_{R,\epsilon} &:= \sigma_R \wedge \sigma_\epsilon, &
  \rho_{R,\epsilon} &:= \rho_R \wedge \rho_\epsilon,&
  \theta_{R,\epsilon} &:= \theta_R \wedge \theta_\epsilon. 
\end{alignat*}

The first assumption is a local Lipschitz condition on the drift
and diffusion co-efficients, together with moment bounds on the
true and numerical solutions.

\begin{assum}
For each $R>0$ there exists a constant $C_R,$ depending only on $R$, such that 
\begin{equation}
|f(a) - f(b)|^2 \vee |g(a) - g(b)|^2 \le C_R |a - b|^2,
    \, \forall a,b \in \RR^m \, \mathrm{with}  \, |a| \vee |b| \le R.
\label{locall}
\end{equation}
For some $p >2$ there is a constant $A$, uniform in $\tau \to 0$, such that
\begin{equation}
\E \left[ \sup_{0 \le t \le T} |\bX(t)|^p \right] \vee
\E \left[ \sup_{0 \le t \le T} |x(t)|^p \right]  \le A.
\label{mbd}
\end{equation}
\label{ass1}
\end{assum}

Note that inequality (\ref{locall}) is a local Lipschitz assumption
which will be satisfied for any $f$ and $g$ in $C^2$.
The inequality (\ref{mbd}) states that the $p^{\rm th}$ moments of
the exact and numerical solution are bounded
for some $p > 2$.  Theorem \ref{t:moment} proves \eqref{mbd}
for the numerical interpolant,  under natural
assumption on $f$ and $g$ (see Assumption \ref{ass:coerce}).
Under the same assumptions, such a bound is known to
hold for $x(t);$ see \cite{mao}.

We clearly also need an assumption on the local error estimate since
if, for example, the drift term $f(u)$ were constant then $E(u,h)\equiv 0$
and the stepsize would, through doubling, reach $\Delta_{\rm max}$, 
no matter how small $\tau$ is, and convergence cannot occur as $\tau \to 0.$
Because the function $F_1(u)$ maps $\RR^m$ into itself, the following
assumption on the zeros of $F_1(u)$ will hold for generic drift functions
$f$ which are non-constant on any open sets; it does exclude, however, the
case of constant drift. Furthermore the assumption on the hitting time
rules out dimension $m=1.$ 
\begin{assum}
\label{ass3}
Define
\begin{equation*}
\ell(\epsilon,R)=d_{H}\{ \Psi(2\epsilon)^c\cap B_R,\Psi(\epsilon)
\cap B_R\} .
\end{equation*} 
For any given $R>0$ we assume that $\ell(\epsilon,R)>0$ for all 
sufficiently small $\epsilon>0$,
and that $\ell(\epsilon,R) \to 0$ as $\epsilon \to 0.$
Furthermore, the hitting time $\rho_{\epsilon}$ satisfies, for any
$X \notin \Psi(0)$,
\begin{equation*}
  \rho_{\epsilon} \to \infty\;{\hbox {as}}\; \epsilon \to 0\;a.s.  
\end{equation*}
\end{assum}
Here $d_{H}$ denotes Hausdorff distance.
The preceding assumption requires that the contours
defining the boundary of $\Psi(\epsilon)$ are strictly nested as $\epsilon$
increases, and bounded.  This enables us to show that the probability of
$(x(t),\bX(t)) \in (\Psi(2\epsilon)^c\cap B_R)
\times (\Psi(\epsilon) \cap B_R)$ is small, a key ingredient
in the proof.

We now state the strong convergence of the adaptive numerical method, using the
continuous-time interpolant $\bX(t).$ 
Note that we do not assume $\Delta_{max} \to 0$ for this theorem.
Hence the non-standard part of the proof comes from estimating
the contribution to the error from regions of phase space where
the time-step is not necessarily small as $\tau \to 0.$

\begin{theorem} Assume that $X \notin \Psi(0).$
Let Assumptions~\ref{ass1} and \ref{ass3} hold. Then, there is $\Delta_c(\tau)$
such that, for all $\Delta_{-1}<\Delta_c(\tau)$ and any $T > 0,$
the numerical solution with continuous-time
extension $\bX(t)$ satisfies
\begin{equation*}
\E \left[ \sup_{0 \le t \le T} |\bX(t) - x(t)|^2 \right]  \to 0
\quad{\hbox {as}}\quad \tau \to 0. 
\end{equation*}
\label{res1}
\end{theorem}

\section{Proof of Convergence Result} \label{sec:alg}

The primary technical difficulty to address in convergence proofs
for adaptive methods is to relate the time-step to the tolerance $\tau$.
Roughly speaking the formula \eqref{eq:form} shows that, 
provided $F_1(u) \ne 0$,
the error control will imply $\Delta_n={\mathcal O}(\tau).$ We now make
this precise.
We provide an upper bound on the timestep sequence of numerical
solutions that remain within $B_{R,\epsilon}$,
for sufficiently small $\tau$.  For given $R, \epsilon >0$ we define the quantities
\begin{equation*} 
  \overline{h}_{R, \epsilon} = \frac{\epsilon}{6K_R} \quad {\rm and}
  \quad  \tau_{R,\epsilon} = \frac{\epsilon^2}{12K_R}.
\end{equation*}

\begin{lemma}\label{lem:ub}
For any $R,\epsilon >0$, if $\{x_n\}_{n=0}^{N}\subseteq
B_{R,\epsilon}$, $\tau < \tau_{R,\epsilon}$ and $\Delta_{-1} <
\frac{2\tau}{\epsilon}$ then
\begin{equation} 
  \Delta_n \leq \min\{\overline{h}_{R,\epsilon}, \frac{2\tau}{\epsilon}\}
  \quad \forall\; 0 \leq n \leq N.\label{le1} 
\end{equation}
\end{lemma}

\begin{proof}
The error control implies
\begin{equation*}
|{E}(x_n, \Delta_n)| = \Delta_n\left|F_1(x_n)+ \Delta_n 
F_2(x_n, \Delta_n)\right| \leq \tau.  
\end{equation*}
Note that 
$$\Delta_{-1}<2\tau_{R,\epsilon}/\epsilon= \overline{h}_{R, \epsilon}.$$
We first proceed by contradiction to prove $ \Delta_n 
\leq\overline{h}_{R,\epsilon}\;\; \forall\; 0 \leq n \leq N $. 
Let $0 \leq m \leq N$ be the first integer 
such that $\Delta_m > \overline{h}_{R,\epsilon}$. 
Then, since there is a maximum timestep
ratio of 2, we have
\begin{eqnarray*} \Delta_m \in (\frac{\epsilon}{6K_R},\frac{\epsilon}{3K_R}]
&  \Rightarrow &
\Delta_m | F_2(x_m,\Delta_m)| < \frac{\epsilon}{2}\\
& \Rightarrow & |{E}(x_m, \Delta_m)| > \Delta_m (\epsilon - 
\epsilon/2) \geq \frac{\epsilon \overline{h}_{R,\epsilon}}{2} = 
\frac{\epsilon^2}{12K_R}=\tau_{R,\epsilon} > \tau. \end{eqnarray*}
Thus $\Delta_m$ is not an acceptable timestep, contradicting our original
assumption. The first result follows. The proof of the bound on the
timestep in (\ref{le1}) now follows immediately since 
\begin{equation*}
\Delta_n \leq \frac{\tau}{|F_1(x_n)+ \Delta_n 
F_2(x_n, \Delta_n)|} \leq \frac{\tau}{(\epsilon-\epsilon/2)} \leq 
\frac{2\tau}{\epsilon}\quad \forall \;0 \leq n \leq N.  
\end{equation*}
\end{proof}

\noindent{\sc Proof of Theorem 1}
We denote the error by
\begin{equation*}
e(t):= \bX(t)-x(t).  
\end{equation*}
Recall the Young inequality: 
for $r^{-1}+q^{-1}=1$
\begin{equation*}
ab \le \frac{\delta}{r}a^r+\frac{1}{q\delta^{q/r}}b^q, \quad
\forall a,b,\delta>0.  
\end{equation*}
We thus have for any $\delta >0$
\begin{eqnarray}
\label{eq:ft0}
\E \left[ \sup_{0 \le t \le T} |e(t)|^2 \right]
&=& \E \left[ \sup_{0 \le t \le T} |e(t)|^2  
            \mathbf{1}\{\theta_{R,\epsilon} > T\}
            \right] + \E \left[ \sup_{0 \le t \le T} |e(t)|^2  
            \mathbf{1}\{\theta_{R,\epsilon} \le T\} \right]
              \nonumber \\
          & \le &
\E \left[ \sup_{0 \le t \le T} |e(t \wedge \theta_{R,\epsilon})|^2
\mathbf{1}\{ \theta_{R,\epsilon} > T\} \right]
       +
     \frac{2 \delta}{p}
     \E \left[ \sup_{0 \le t \le T} |e(t)|^p \right] \nonumber \\
   && \mbox{} + \frac{1-\FFrac{2}{p}}{\delta^{2/(p-2)}}
             \Prob \big( \theta_{R,\epsilon} \le T \big). 
\end{eqnarray}
Now 
$$\Prob \left( \theta_{R,\epsilon} \le T \right) 
= \Prob\{\theta_R \le T\}+\Prob\{\theta_{\epsilon}\le T, \theta_R>T\}.
$$
But
\begin{equation*}
  \Prob \{\theta_{R}\le T\} \le \Prob\{\sigma_R \le T\}+\Prob\{\rho_R \le T\}
\end{equation*}
whilst
\begin{equation*}
  \Prob\{\theta_{\epsilon}\le T, \theta_R>T\}  \le
  \Prob\{\rho_{\epsilon} \le T\}+
  \Prob\{\theta_{\epsilon}\le T, \theta_R>T,\rho_{\epsilon} > T\}.
\end{equation*}
Thus we have
\begin{equation*}
\Prob \left( \theta_{R,\epsilon} \le T \right) 
\le \Prob(\sigma_R \le T)+\Prob (\rho_R \le T)+\Prob(\rho_{\epsilon} \le T)+
\Prob\{\theta_{\epsilon}\le T, \theta_R>T,\rho_{\epsilon} > T\}.
\end{equation*}
To control the last term notice that whenever $\theta_{\epsilon}\le
T$,  $\theta_R>T$ and $\rho_{\epsilon} > T$ we know that
$|e(\sigma_{\epsilon})| \ge \ell(\epsilon,R)$. Hence we have
\begin{equation*}
\Prob\{\theta_{\epsilon}\le T, \theta_R>T,\rho_{\epsilon} > T\} \le 
\Prob \{|e(T \wedge \theta_{R,\epsilon})| \ge \ell(\epsilon,R)\} \le
\E |e(T \wedge \theta_{R,\epsilon})|^2/\ell(\epsilon,R)^2.
\end{equation*}
Combining the two preceding inequalities gives
\begin{equation*}
\Prob \left( \theta_{R,\epsilon} \le T \right) 
\le \Prob(\sigma_R \le T)+\Prob (\rho_R \le T)+\Prob(\rho_{\epsilon} \le T)+
\E \Bigr(\sup_{0 \le t \le T}|e(t \wedge \theta_{R,\epsilon})|^2\Bigl)/\ell(\epsilon,R)^2.
\end{equation*}
By Markov's inequality
\begin{equation*}
\Prob\{\sigma_R \le T\}, \Prob\{\rho_R \le T\} \le \frac{A}{R^p}.  
\end{equation*}
so that
\begin{equation}
\Prob \left( \theta_{R,\epsilon} \le T \right) 
\le \frac{2A}{R^p}+\Prob(\rho_{\epsilon} \le T)+
\E \Bigr(\sup_{0 \le T}|e(t \wedge \theta_{R,\epsilon})|^2\Bigl)/\ell(\epsilon,R)^2.
\label{eq:now1}
\end{equation}
Furthermore,
\begin{equation}
\label{eq:now2}
  \E \left[ \sup_{0 \le t \le T} | e(t)|^p \right]
     \le {2}^{p-1}
   \E \left[ \sup_{0 \le t \le T} \left(|\bX(t)|^p + |x(t)|^p \right) \right]
    \le 2^p A.
\end{equation}

Using \eqref{eq:now1}, \eqref{eq:now2}
in (\ref{eq:ft0}) gives, for $\epsilon$ sufficiently small,
\begin{eqnarray}
\E \left[ \sup_{0 \le t \le T} |e(t)|^2 \right] &\le&
\left(1+\frac{p-2}{p\delta^{2/(p-2)}\ell(\epsilon,R)^2}\right)
\E \left[ \sup_{0 \le t \le T} |e(t \wedge \theta_{R,\epsilon})|^2 \right]
   \nonumber \\
       && \mbox{}
     + \frac{2^{p+1} \delta A}{p}  +
           \frac{(p-2)}{p\delta^{2/(p-2)}}\left[\frac{2A}{R^p}+
\Prob\{\rho_\epsilon \le T\}\right].
   \label{eq:bdd}
\end{eqnarray}

Take any $\kappa>0$.
To complete the proof we choose $\delta$ sufficiently small
so that the second term on the right hand side of
(\ref{eq:bdd}) is bounded by $\kappa/4$
and then $R$ and $\epsilon$ sufficiently large/small so that the
third and fourth terms are bounded by  $\kappa/4.$ Now reduce $\tau$ 
so that Lemma \ref{lem:ub} applies. Then, by further reduction of $\tau$ 
in Lemma \ref{lem:basic}, we upper-bound the first term by $\kappa/4.$ 
(Lemma \ref{lem:basic} calculates the error conditioned
on the true and numerical solutions staying within a ball of radius $R$,
and away from small sets where the error control mechanism breaks down.
With this conditioning it follows from Lemma \ref{lem:ub} that
we have $\Dt_n={\calO}(\tau)$, which is the essence of why
Lemma \ref{lem:basic} holds.) 

Consequently we have
\begin{equation*}
  \E \left[ \sup_{0 \le t \le T} | \bX(t) - x(t) |^2 \right] \le \kappa
\end{equation*}
and since $\kappa$ is arbitrary the required result follows.
\epf

In the following, $C$ is a universal
constant independent of $T,R$, $\epsilon, \delta$ and $\tau$. 
Likewise $C_R$ is a universal constant depending upon $R$, but independent
of $T,\epsilon,\delta$ and $\tau,$ 
$C_{R,T}$ is a universal constant depending upon $R$ and
$T$, but independent of $\epsilon, \delta$ and $\tau$ and $C_{R,\epsilon,T}$
and so forth are defined similarly. 
The actual values of these constants may change from one occurrence to the 
next.

\begin{lemma}
\label{lem:basic}
Assume that $X \notin \Psi(0)$ and that $\tau$ is sufficiently small for
the conditions of Lemma
\ref{lem:ub} to hold. Then the continuous interpolant of the numerical
method, $\bX(t)$, satisfies the following error bound:
\[ \E \left[ \sup_{0 \le t \le T} |\bX(t \wedge \theta_{R,\epsilon})
 - x(t \wedge \theta_{R,\epsilon})|^2 \right] \le  C_{R,\epsilon,T}\tau.\]
\end{lemma}

\begin{proof}
Using
\begin{equation*}
  x(t \wedge \theta_{R,\epsilon}) := X +
\int_{0}^{t \wedge \theta_{R,\epsilon}} f(x(s)) ds +
\int_{0}^{t \wedge \theta_{R,\epsilon}} g(x(s)) dW(s),
\end{equation*}
equation (\ref{eq:bx1}) and Cauchy--Schwartz, we have
that $\chi:=| \bX(t \wedge \theta_{R,\epsilon}) - x(t \wedge
\theta_{R,\epsilon}) |^2$, satisfies
\begin{eqnarray*}
   \chi&=& \left| \int_{0}^{t \wedge \theta_{R,\epsilon}}
\Bigl(f(X(s)) - f(x(s))\Bigr) ds  + \int_{0}^{t \wedge \theta_{R,\epsilon} }
\Bigl( g(X(s)) - g(x(s))\Bigr) dW(s)\right|^2 \\ & \le & 2 \left[T
                \int_{0}^{t \wedge \theta_{R,\epsilon}}
| f(X(s)) - f(x(s)) |^2 ds 
+ \left| \int_{0}^{t \wedge \theta_{R,\epsilon}} \Bigl(g(X(s)) -
 g(x(s))\Bigr) dW(s) \right|^2 \right].
\end{eqnarray*}  
Let
$$E(s):=
\left[ \sup_{0 \le t \le s} | \bX(t \wedge \theta_{R,\epsilon}) -
x(t \wedge \theta_{R,\epsilon}) |^2 \right]
$$
Then, from the local Lipschitz condition (\ref{locall}) and
the   Doob-Kolmogorov Martingale inequality \cite{RW}, we have for any $t^* \le T$
\begin{eqnarray}  
\E E(t^*) 
&\le& 2C_R(T + 4) \E
\int_{0}^{t^* \wedge \theta_{R,\epsilon}} | X(s) - x(s) |^2 ds
\nonumber \\ & \le & 4 C_R(T + 4) \E
\int_{0}^{t^* \wedge \theta_{R,\epsilon}}   
\Bigl[ |X(s) - \bX(s)|^2  + | \bX(s) - x(s) |^2 \Bigr] ds
\nonumber \\
& \le &
4 C_R(T + 4) \left[ \E \int_{0}^{t^* \wedge \theta_{R,\epsilon}}
 |X(s) - \bX(s)|^2  ds +\int_0^{t^*} \E E(s)ds \right].  
     \label{bdi}
\end{eqnarray}
Given $s \in [0,T \wedge \theta_{R,\epsilon})$,
let $k_s$ be the integer for which $s \in [t_{k_s},t_{k_s+1})$.
Notice that $t_{k_s}$ is a stopping time
because $\Delta_{k_s}$ is a deterministic function of 
$(x_{k_s},\Delta_{k_s-1})$.  
We now bound the right hand side in \eqref{bdi}.
From the local Lipschitz condition (\ref{locall}),
a straightforward calculation shows that
\begin{equation*}
  | X(s) - \bX(s) |^2 
  \le  C_R(|x_{k_s}|^2 + 1)( \Delta_{k_s}^2 +
  | W(s) - W(t_{k_s}) |^2 ). 
\end{equation*}
Now, for $s<\theta_{R,\epsilon}$, using Lemma \ref{lem:ub},
\begin{align*}
|W(s)-W(t_{k_s})|^2 &= s-t_{k_s}+2\int_{t_{k_s}}^s
[W(l)-W(t_{k_s})]dW(l)\\
& \le (s-t_{k_s})[1+I(s)] \le \frac{2\tau}{\epsilon}[1+I(s)].
\end{align*}
Here
$$I(s)=\frac{2}{(s-t_{k_s})}\Big|\int_{t_{k_s}}^s
[W(l)-W(t_{k_s})]dW(l)\Big|.$$
Let ${\mathcal H}_s$ denote the $\sigma-$algebra of Brownian paths up to time
$t_{k_s}$. Then, conditioned on ${\mathcal H}_s$, we have
\begin{equation}
\label{eq:saywhat?}
\E I(s) \le \sqrt 2.
\end{equation}
Thus, using Lemma \ref{lem:ub}, (\ref{mbd}) and the Lyapunov inequality 
\cite{KP},
\begin{eqnarray*}
  \E \int_{0}^{t^* \wedge \theta_{R,\epsilon}}  | X(s) - \bX(s) |^2 ds
  &\le& \E \int_{0}^{t^* \wedge \theta_{R,\epsilon}  }
   C_R (|x_{k_s}|^2 + 1 )( 4\tau^2/\epsilon^2 +
  | W(s) - W(t_{k_s}) |^2 ) ds \\
&\le & C_{R,\epsilon}\tau \E \int_{0}^{t^*}
(1+|x_{k_s}|^2)(1+I(s))ds\\
&\le& C_{R,\epsilon,T}(A^{2/p+1})\tau.
\end{eqnarray*}
To obtain the last line we condition on ${\mathcal H}_s$ so that
$|x_{k_s}|^2$ and $I(s)$ are independent; we then use \eqref{eq:saywhat?}
and the assumed moment bound.

In (\ref{bdi}), we then have by Lemma \ref{lem:ub}
\begin{eqnarray*}
\E E(t^*) 
&\le&  C_{R,\epsilon, T} \tau+4 C_{R,T} \int_{0}^{t^*} \E E(s)ds. 
\end{eqnarray*}
Applying the Gronwall inequality we obtain
\begin{equation*}
\E \left[ \sup_{0 \le t \le T} \left( \bX(t \wedge \theta_{R,\epsilon})
 - x(t \wedge \theta_{R,\epsilon}) \right)^2 \right]
\le  C(R,\epsilon,T)\tau.  
\end{equation*}
\end{proof}

\section{Stability Results}\label{sec:mainResults}

For all of our stability results, in this and the following
sections, we make the assumption that
\eqref{eq:coerce} holds, together with some conditions on the diffusion
matrix. To be explicit we make 
\begin{assum}
\label{ass:coerce}
There exists finite positive $\alpha, \beta$ such that
\begin{equation*}
\la f(x),x \ra \le \alpha-\beta |x|^2 \quad \forall x \in \real^m,
\end{equation*}
where $\la \cdot, \cdot \ra$ is the inner-product inducing the 
Euclidean norm $|\cdot|$.
Furthermore $m=d$, $\Sigma$ is globally bounded  and is globally invertible.
\end{assum}

The assumption is made, without explicit statement, for the remainder of the paper. We
also assume, without explicit statement, that $\tau<2\beta$ so that $\tb>0.$
Finally we assume, 
also without explicit statement, that there is at least one point
$\by \in \real^m$ such that
\begin{equation}
\label{eq:ratz2}
k^*(\by)={\mathcal G}(\by,1).
\end{equation}
This may implicitly force upper bounds on $\tau$ and $\Delta_{max}$, although
neither is necessarily restricted by this assumption. The existence of such a $\by$
is implied by Assumption \ref{ass:coerce}, which rules out $f$ being identically
constant. Then there exists $\by$ for which the function
$$|f(\by+hf(\by))-f(\by)|$$
is non-zero in a neighbourhood of $h=0$ and \eqref{eq:ratz2} must hold,
possibly after enforcing bounds on $\tau$ and $\Delta_{max}.$

Under Assumption \ref{ass:coerce} the solution
of \eqref{eq:sde} exists for all $t>0$ \cite{has, mao} and
the equation is geometrically ergodic \cite{has, MT, msh}.
The first stability result ensures that the method will not
decrease its stepsize in such a way that it is unable to reach
arbitrary finite times.
\begin{theorem}
\label{t:stf}
The stopping times $N_j$ are almost surely finite.
\end{theorem}
The next result is the main ergodic result of the paper. It
ensures that the adaptive method has an attracting statistical steady
state.  Letting $\EE^{y,l}$ denote the expectation under the
Markov chain started at $x_0=y, k_{-1}=l$, we have the following result.
(Recall $\delta$ occurring in the definition of stopping times $N_j.$)

\begin{theorem} \label{t:erg}  Assume that $\delta>5\Delta_{max}.$ 
The Markov chain $\{y_j,l_j\}=\{x_{N_j+1},k_{N_j}\}$ has a unique
invariant measure $\pi$. Furthermore, if $h:\real^m \times \ZZ^+ \to \real$
is measurable and
\begin{equation*}
  |h(y,l)| \le 1+|y|^2 \quad \forall (y,l) \in \real^m \times \ZZ^+,
\end{equation*}
then there exists $\lambda \in (0,1)$, $\kappa \in (0,\infty)$
such that
\begin{equation*}
  |\EE^{y_0,l_0} h(y_n,l_n)-\pi(h)| \le \kappa \lambda^n[1+|y_0|^2].  
\end{equation*}
\end{theorem}

The final result gives a moment bound on the continuous time interpolants
of the numerical solution, mimicking that for the SDE itself.
\begin{theorem}
\label{t:moment} There exists a $\lambda >0$ and a $c>0$ so that
\begin{align*}
\bbE \exp( \lambda \sup_{t \in [0,T]} \|X(t)\|^{2}) &\le \exp( \lambda
|X|^2 + c T)\\
\bbE \exp( \lambda \sup_{t \in [0,T]} \|\bx(t)\|^{2}) &\le \exp( \lambda
|X|^2 + c T).
\end{align*}
\end{theorem}

\section{Proof of Stability Results}\label{sec:APriori}
\setcounter{equation}{0}

We start with a number of estimates which will be needed to prove the
main results. It is useful to define
\begin{align*}
\xi_{n+1}=2\sqrt{\Dt_{n}}\la \xs_n,\Sigma(x_n)\eta_{n+1} \ra,&\quad
\xit_{n+1}=\Dt_{n}[|\Sigma(x_n)\eta_{n+1}|^2-\sigma^2(x_n)],\\
M_n=\sum_{j=0}^{n-1}\xi_{j+1},&\quad \Mt_n=\sum_{j=0}^{n-1}\xit_{j+1}.
\end{align*}
Observe that $\la \xs_n,\Sigma(x_n)\eta_{n+1} \ra$ is a Gaussian
random variable conditioned on the values of $x_n$ and  $\xs_n$. Hence
the last two
expressions are Martingales satisfying the assumptions of Lemma
\ref{l:expMartNew} from the appendix. Also notice that the quadratic variations satisfy
\begin{align}
\label{eq:oreza}
\la M \ra_n \le    \sum_{j=0}^{n-1}
4\Dt_{j}|\xs_j|^2\sigma^2 \quad\mbox{and}\quad
\la \Mt \ra_n \le   \sum_{j=0}^{n-1} a\Dt_{j}^2\sigma^4.
\end{align}

We start with a straightforward lemma.
\begin{lemma} The sequences $\{\xs_n\}$ and $\{x_n\}$ satisfy 
  \label{l:4.1}
  \begin{align*}
    |\xs_{n}|^2 &\le |x_n|^2+2\Dt_{n}[\ta-\tb |\xs_n|^2],\\
    |x_{n+1}|^2 &\le \beta_n|x_n|^2+\Dt_{n}[2\ta+\sigma^2]
+\xi_{n+1}+\xit_{n+1}.
  \end{align*}
Hence 
  \begin{equation*}
    \la M \ra_n \le 4 \sigma^2\sum_{j=0}^{n-1} |x_j|^2\Dt_{j}+
   8\sigma^2\ta \sum_{j=0}^{n-1} \Dt_{j}^2.  
  \end{equation*}
\end{lemma}
\vspace{0.1in}

\bpf Taking the inner product of the equation
\begin{equation*}
  \xs_n=x_n+\Dt_{n} f(x_n)  
\end{equation*}
with $\xs_n$ and using the fact that
the error control implies
\begin{equation*}
  |f(\xs_n)-f(x_n)| \le \tau,  
\end{equation*}
a straightforward calculation
from \cite{SH2}, using \eqref{eq:coerce}, gives the first result.
To get the second simply square the expression \eqref{eq:ratz}
for $x_{n+1}$ and use the first, noting that $\beta_n \le 1.$
For the third use the first in the bound \eqref{eq:oreza} for $\la M \ra_n$. 
\epf

\begin{lemma}
  \label{l:growSlow}
  We have
  \begin{align*}
    |x_{n+1}|^2 \leq |X|^2 + C_0 t_{n+1} + M_{n+1} -\frac12 \frac{\tilde
      \beta}{\sigma^2} \langle M \rangle_{n+1} + \tilde M_{n+1} -
    2\langle \tilde M \rangle_{n+1}  \; .
  \end{align*}
  where $C_0=[2\tilde \alpha + 4\sigma^4\Delta_{max}]$. Furthermore
  \begin{align*}
    \PP\Bigl( \sup_{0\leq n} \{|x_n|^2 - C_0 t_n\} \geq |X|^2 +
    A \Bigr) \leq 2 \exp\big( -B A \big)
  \end{align*}
  where $B$ is a positive constant depending only on $\sigma$
  and $\tilde \beta$. 
\end{lemma}
\bpf
Squaring the expression for $x_{n+1}$ in \eqref{eq:ratz},
bounding $|\xs_n|^2$ by the first inequality in Lemma  \ref{l:4.1} and summing gives
\begin{align*}
  |x_{n+1}|^2 \leq |X|^2 +C_0 t_{n+1} + S_{n+1} +
   \tilde S_{n+1}
\end{align*}
where
\begin{equation*}
S_{n+1}=M_{n+1} - 2 \tilde \beta \sum_{k=0}^n |\xs_k|^2 \Delta_{k}, \quad
  \tilde S_{n+1}=\tilde M_{n+1} - 4 \sigma^4 \Delta_{max} t_{n+1}  
\end{equation*}
and $M_n$, $\tilde M_n$ are as before. Using \eqref{eq:oreza}, one obtains
\begin{align*}
\langle \tilde M \rangle_{n+1} \leq 2\sigma^4 \Delta_{max} t_{n+1}\;.
\end{align*}
and 
\begin{align*}
  \langle M \rangle_{n+1} 
  \leq 4 \sigma^2 \sum_{k=0}^n \Delta_{k} |\xs_k|^2\; . 
\end{align*}

Combining all of this produces,
\begin{align*}
  S_{n+1} \leq M_{n+1} -\frac12 \frac{\tilde \beta}{\sigma^2} \langle M
\rangle_{n+1} \quad\mbox{ and }\quad 
\tilde S_{n+1} \leq \tilde M_{n+1} - 2\langle \tilde M \rangle_{n+1} \; .  
\end{align*}
The probabilistic estimate follows from the exponential martingale
estimates from the Appendix.
\epf 

\begin{corollary}\label{c:growSlow}
  Then there exists a universal $\lambda >0$ and $C_1>0$ so that for any stopping
  time $N$ with $0\leq t_N\leq t_*$ almost surely, for some fixed
  number $t_*$, one has 
  \begin{align*}
    \EE \exp( \lambda \sup_{0 \leq n \leq N} |x_n|^2 ) \leq
    C_1\exp(\lambda |X|^2 + \lambda C_0t_*)\;.
  \end{align*}
\end{corollary}
\begin{proof} The result follows from Lemma \ref{l:growSlow} and the
  observation that 
  \begin{equation*}
    \PP\Bigl( \sup_{0\leq n \leq N} |x_n|^2 \geq |X|^2 +C_0 t_*+A
    \Bigr)\leq  \PP\Bigl( \sup_{0\leq n} \{|x_n|^2 - C_0 t_n\} \geq
    |X|^2 +     A \Bigr)\;.
\end{equation*}
\end{proof}

\vspace{0.1in}

\begin{lemma} The Markov chain $\{x_{N_j}\}_{j \in \ZZ^+}$ satisfies the
Foster-Lyapunov drift condition
\begin{equation*}
  \bbE \{|x_{N_{j+1}}|^2 | \cF_{N_j}\} \le
  \exp(-2\gamma^{-} \tb \dm) |x_{N_j}|^2+\exp(2\tb \delta^+)[2\ta+\sigma^2]\dpp.  
\end{equation*}
That is
\begin{equation*}
  \bbE \{|y_{j+1}|^2 | \cG_{j}\} \le
  \exp(-2\gamma^- \tb \dm) |y_{j}|^2+\exp(2\tb \delta^+)[2\ta+\sigma^2]\dpp.  
\end{equation*}
\label{l:4.2} 
\end{lemma}

\vspace{0.1in}

\bpf Note that $(1+x)^{-1} \le e^{-\gm x}$ for all $x \in
[0,2\Delta_{max}\tb]$. From Lemma \ref{l:4.1}, we have
$$|x_{n+1}|^2 \le \beta_n |x_n|^2+\kk_n+\xi_{n+1}+\xit_{n+1}$$
where $\kk_n:=\Dt_{n}[2\ta+\sigma^2]$. Defining
$$\gamma_j=\Bigl(\Pi_{l=0}^{j-1} \beta_l^{-1} \Bigr)$$
we obtain 
\begin{equation*}
\bbE \bigl(\gamma_{N_{j+1}}  |x_{N_j+1}|^2|\cF_{N_j}\bigr)  \le 
\gamma_{N_j}  |x_{N_j}|^2 + \bbE \bigl(\sum_{l=N_j}^{N_{j+1}-1} \gamma_{l+1}
\kk_l|\cF_{N_j}\bigr). 
\end{equation*}
Now 
\begin{equation}
  \sum_{l=N_j}^{N_{j+1}-1} \Dt_{j}  \le  \delta+\Delta_{max}  = \dpp
  \qquad\mbox{and}\quad \sum_{l=N_j}^{N_{j+1}-1} \Dt_{j} \ge  \delta  = \dm.  
\label{eq:enothis}
\end{equation}
Straightforward calculations show that
$$\gamma_{N_{j+1}} \ge \exp(2\tb \gamma^- \delta) \gamma_{N_j}$$
and
$$\gamma_{l+1} \le \exp(2\tb \delta^+) \gamma_{N_j}.$$
Hence
\begin{equation*}
\bbE \big\{ |x_{N_j+1}|^2| \cF_{N_j}\big\} \le 
\exp(-2 \gamma^- \tb \dm)|x_{N_j}|^2
+\exp(2\tb \delta^+) \bbE\Big\{ \sum_{l=N_j}^{N_{j+1}-1} \kk_l | 
\cF_{N_j} \Big\}  \end{equation*}
and for the required result we need to bound the last term. By
\eqref{eq:enothis} we have
\begin{equation*}
\bbE \Bigl(\sum_{l=N_j}^{N_{j+1}-1} \Dt_{l}|\cF_{N_j}\Bigr) \le \delta+\Delta_{max}=\dpp  
\end{equation*}
and we obtain
\begin{equation*}
  \EE \Bigl(\sum_{l=N_j}^{N_{j+1}-1} \kk_l |\cF_{N_j}\Bigr) \le [2\ta+\sigma^2]\delta^{+}.
\end{equation*}
This gives the desired bound.
\epf

\vspace{0.1in}

We now proceed to prove the ergodicity and moment bound.
We prove geometric ergodicity of the Markov chain $\{y_j,l_j\}$
by using the approach highlighted in \cite{MT}. In particular
we use a slight modification of Theorem 2.5 in \cite{msh}. Inspection
of the proof in the Appendix of that paper shows that, provided an invariant
probability measure exists, and this follows from Lemma \ref{l:4.2},
the set $C$ in the minorization condition need not be compact: it simply needs
to be a set to which return times have exponentially decaying tails.

Let
\begin{equation*}
  P(y,l,A)=\PP((y_1,l_1) \in A|(y_0,l_0)=(y,l))  
\end{equation*}
where 
\begin{equation*}
  A \in \cB(\real^m) \otimes \cB(\ZZ^+), (y,l) \in \real^m \times \ZZ^+.  
\end{equation*}
We write $A=(A_y,A_l)$ with $A_y \in \cB(\real^m)$ and $A_l \in \cB(\ZZ^+).$

The minorization condition that we use, generalizing that in Lemma 2.5 of
\cite{msh}, is now proved:\footnote{Note that although $C$ is compact
in the following, $C \times \ZZ^+$ is not.}

\begin{lemma} \label{l:minor}
Let $C$ be compact. For $\delta>5\Delta_{max}$ there is $\zeta>0$, $\by \in \real^m$
and $\ddelta>0$ such that 
\begin{equation*}
  P(y,l,A) \ge \zeta \nu(A)\quad \forall A \in \cB(\real^m) \times
  \cB(\ZZ^+),\;(y,l) \in C \times \ZZ^+,
\end{equation*}
where  
\begin{equation*}
  \nu(A)=Leb\Big\{B(\by,\ddelta) \cap A_y\Big\} \cdot 
  \ONE\Big\{\ks(\by) \in A_l\Big\}.  
\end{equation*}
\end{lemma}

\vspace{0.1in}

\bpf Let $M=N_1(2\Delta_{max})$ and $N=N_1(\delta).$ Recall the
definition \eqref{eq:ratz2} of $\bar{y}.$ Since $t_M \leq 3
\Delta_{max}$ almost surely, setting $r^2=R^2+B^2$ Corollary \ref{c:growSlow} implies we
can choose positive $B$ and $R$ sufficiently large so that
\begin{equation*}
  \PP\Big\{\sup_{0 \le n \le M} |x_n|^2 \le r^2\Big\} \ge \frac12  
\end{equation*}
and $\by \in B(0,r)$,  and $C \subseteq B(0,R)$. 
Label this event, with probability in excess of $\frac12$, by $E_1.$
If $E_1$ occurs then there exists $l \in \{0, \dots, M\}$ such that
$k_l \le k^*(B(0,r))$. This follows by contradiction, since otherwise
$\Delta_j=2^{j+1}\Delta_{-1}$ for $j \in \{0,\dots,M\}$ and 
$$t_M = \sum_{j=0}^{M-1} \Delta_j \le \Delta_{max}\sum_{j=0}^{M-1}2^{j-M} 
\le \Delta_{max}\sum_{k=1}^{\infty} 2^{-k}=\Delta_{max}.$$
However $t_M \ge 2\Delta_{max}$, a contradiction.
Once $k_j \le k^*(B(0,r))$ it follows that
$k_n \le k^*(B(0,r))$ for $n \in \{l,\dots,M\}$ as a consequence
of the step-size selection mechanism.

Assume that $E_1$ has occurred.
By choice of $\epsilon$ sufficiently small, $B(\by,\epsilon) \subseteq B(0,r)$.
We now choose the $\eta_j$ for $j \in \{M,\dots,N-1\}$ to
ensure the event $E_2$ namely that
\begin{equation*}
  x_j \in B(\by,\ddelta), \quad M+1 \le j \le N.  
\end{equation*}
It is possible to ensure that the event has
probability $p_1>0$, uniformly for  $X \in C$ and $k_0 \in \ZZ^+.$
The fact that $x_M \in B(0,r)$ gives uniformity in $X \in C.$
We prove an upper bound on the number of steps
after $M$ to get probability independent of $k_{-1} \in \ZZ^+.$
To do this notice that $k_n \le k^*(B(0,r))$ now for
$n \in \{j, \dots,N\}$, again as a consequence of the step-size mechanism.
In fact $k_N=k^*(B(\by,\epsilon))=k^*(\by).$
This follows because an argument analogous to that
above proves that there is $l\in\{M+1,\dots,N\}$ for which
$k_n \le k^*(B(\by,\epsilon))=k^*(\by)$ for $n \in \{l,\dots,N\}.$
Now
$k^*(B(\by,\epsilon))=k^*(\by)$, by continuity of $f$ and possibly by further
reduction of $\epsilon.$ Since $k_j<k^*(\by)={\mathcal G}(\by,1)$ is not possible,
it follows that $k_N=k^*(\by).$ 

If $E_1$ and $E_2$ both occur then,
for some $\gamma>0$, the probability that $y_1=x_{N_1(\delta)} \in A_y$
is bounded below by $\gamma Leb\{A_y \cap B(\by,\ddelta)\}$,
for some $\gamma>0$, because
\begin{equation*}
  x_{N}=\xs_{N-1}+\sqrt{\Dt_{N-1}}
  \Sigma(\xs_{N-1})\eta_{N},  
\end{equation*}
$x_{N-1}$ is in a compact set and $\Sigma$ is invertible.
The fact that $\eta_j$ are i.i.d 
Gaussian gives the  required lower bound in terms
of Lebesgue measure. The final result follows with $\zeta=\gamma p_1/2.$
\epf

\vspace{0.1in}

With this minorization condition in hand, we turn to the proof of
ergodicity.

\vspace{0.1in}

\bpf[Proof of Theorem \ref{t:erg}] The existence of an invariant
measure $\pi$ follows from the Foster-Lyapunov drift condition of
Lemma \ref{l:4.2} which gives tightness of the Krylov--Bogoljubov
measures. Lemma \ref{l:4.2} shows that the chain $\{y_j,l_j\}$
repeatedly returns to the set $C \times \ZZ^+$ and that the return
times have exponentially decaying tails. This generalizes Assumption
2.2 in \cite{msh}.  Lemma \ref{l:minor} gives a minorization condition
enabling a coupling. Together these two results give Theorem
\ref{t:erg}, by applying a straightforward modification of the
arguments in Appendix A of \cite{msh}.  \epf

\vspace{0.1in}

\bpf[Proof Theorem \ref{t:moment}] We define the stopping time $N$ by
$$N=\inf_{n \ge 0}\{n: t_n \ge T\}$$
noting that 
$$T \le t_{N} \le T+\Delta_{max}$$ 
Notice that \eqref{eq:bx2} implies that
$$\sup_{0 \le t \le t_N}|\bx(t)|^{2} \le C\left[
1+\sup_{0 \le k \le N}|x_k|^{2}+\sup_{(t-s) \in [0,\Delta_{max}], s \in [0,T]}
|W(t)-W(s)|^{2}\right].$$
Here we have used the fact that, by Lemma \ref{l:4.1},
$$|\xs_k|^2 \le |x_k|^2+2\Delta_{max}\ta.$$
From this relationship between the supremum of moments of $\bx(t)$ and
$X(t)$, and from the properties of increments of Brownian motion,
it follows that, to prove Theorem \ref{t:moment}, it suffices
to bound 
$$\bbE \exp( \lambda \sup_{0 \le k \le N}|x_k|^{2})$$
for some $\lambda >0$. However this follows from the fact that $t_N
\leq T+\Delta_{max}$ and Corollary \ref{c:growSlow}.

\section{Numerical Experiments: Pathwise Convergence}
\label{sec:numerics}
\setcounter{equation}{0}

We now provide some numerical experiments to complement the analysis
of the previous sections. We begin, in this section,
by demonstrating the importance of Assumption \ref{ass3} in 
ensuring pathwise convergence. In the next
section we discuss an abstraction of the method
presented and studied in detail in this paper. 
Section \ref{sec:9} then shows how this abstraction leads 
to a variant of the method discussed here, tailored to
the study of damped-driven Hamiltonian systems. We
provide numerical experiments showing the efficiency of the
methods at capturing the system's invariant measure.

In the convergence analysis, we made Assumption
\ref{ass3} the second part of which was to assume that the hitting
time of small neighbourhoods of the set $\Psi(0)$ is large with high
probability. We now illustrate that this is not simply a technical
assumption. We study the test problem
\begin{equation}
  \label{eq:test2}
  dy=y-y^3+dW,
\end{equation}
where $W$ is a real valued scalar Brownian motion. 
The set $\Psi(0)$ comprises the points where $f(y):=y-y^3$
satisfies $f(y)=0$ and $f'(y)=0$, that is the points
$\pm 1,0, \pm \frac{1}{\sqrt 3}.$ Since the problem is one dimensional
the hitting time to neighbourhoods of these points is not small.

\begin{figure}[tbhp]
   \begin{tabular}{cc}
    \includegraphics[width=2.75in]{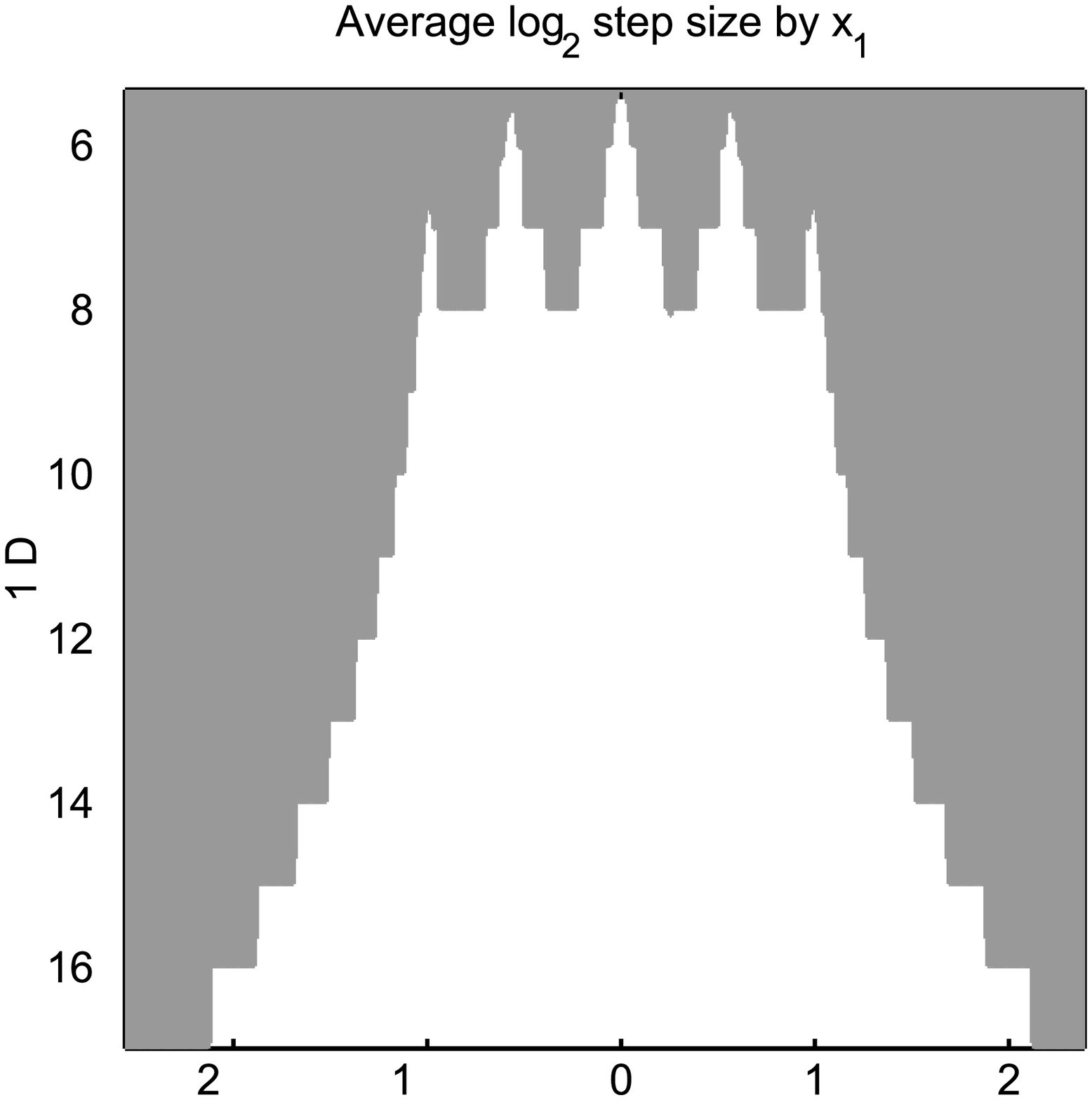} & \includegraphics[width=2.75in]{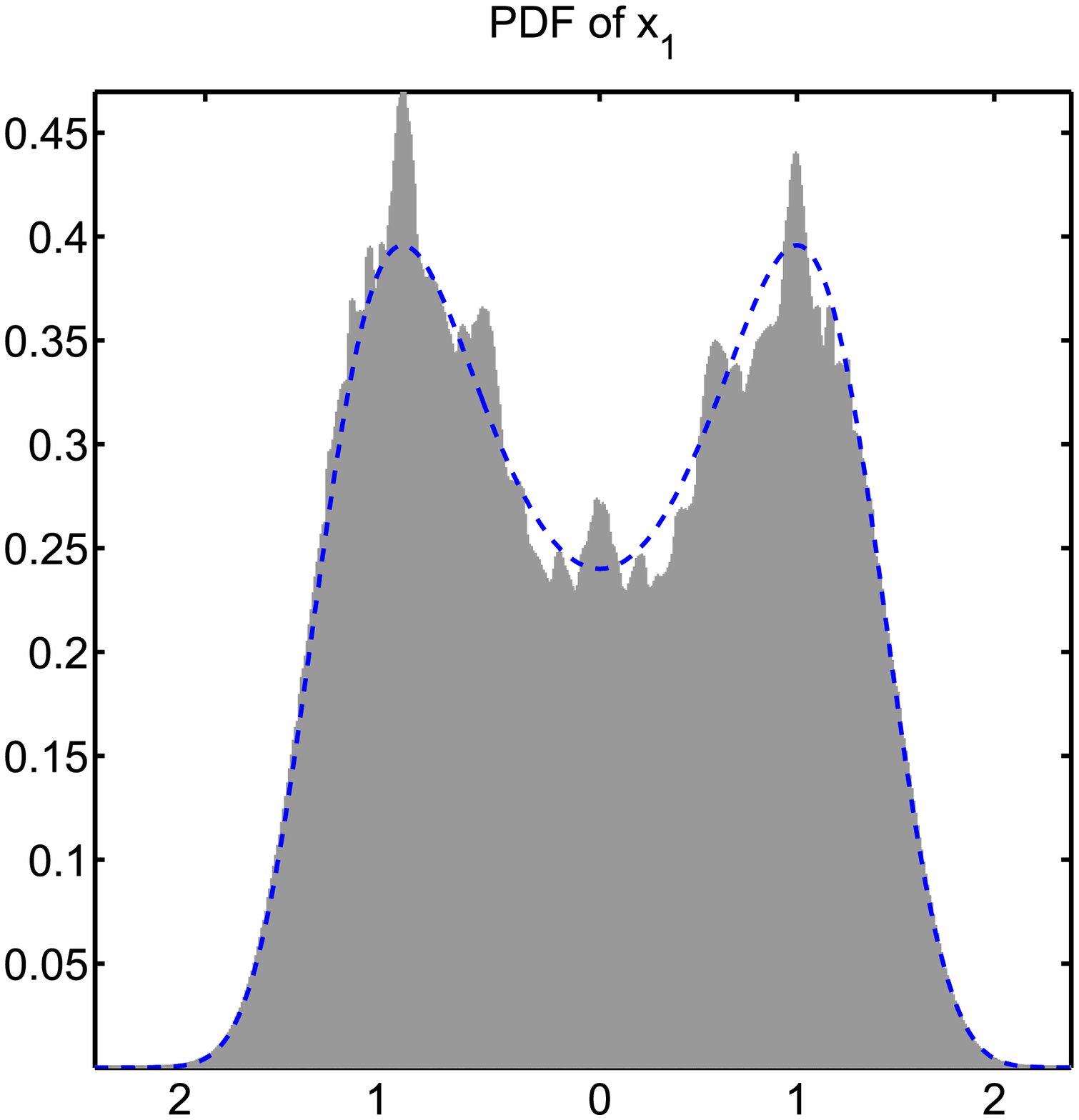}\\
    \includegraphics[width=2.75in]{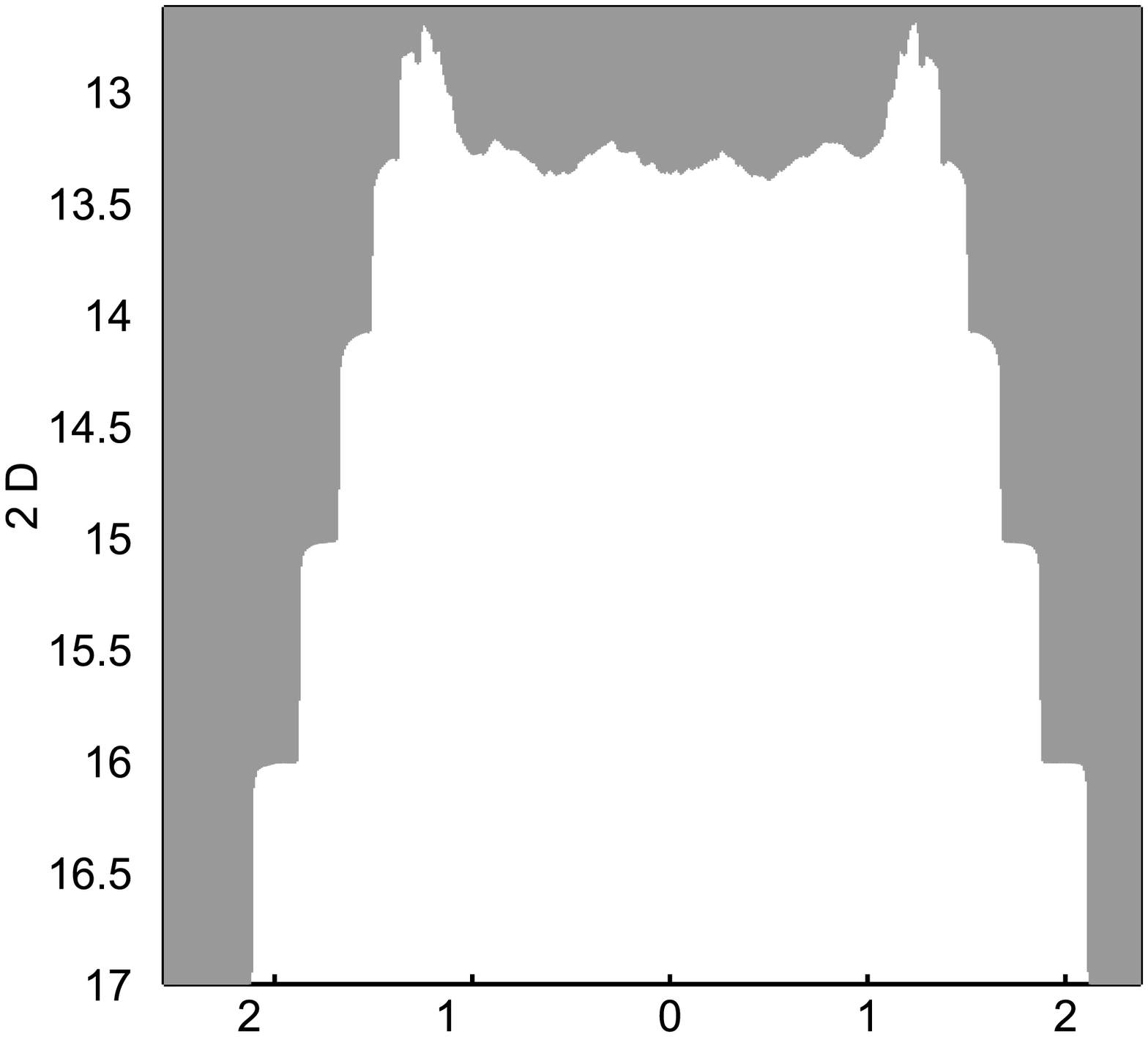} &  \includegraphics[width=2.75in]{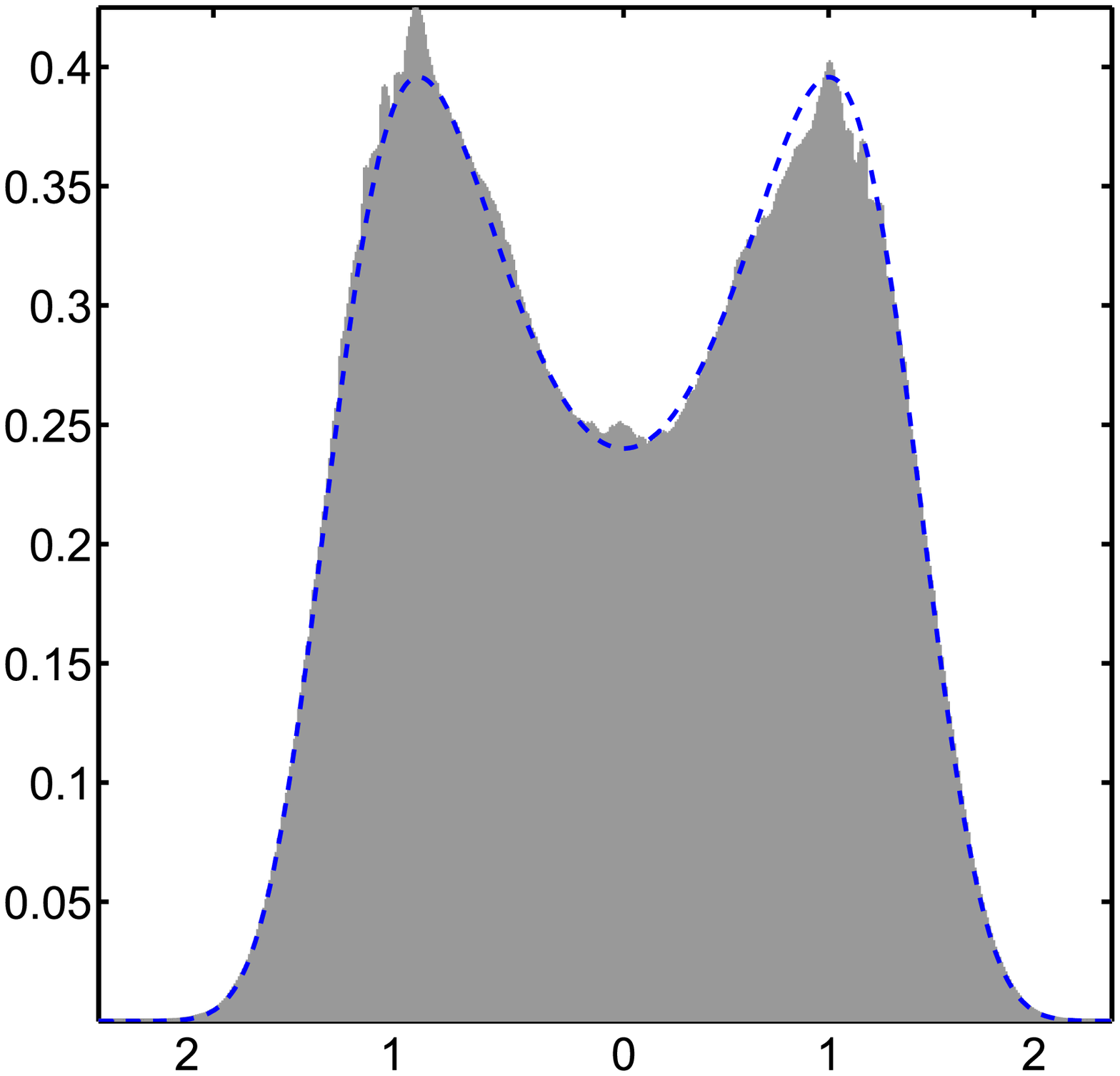}\\
    \includegraphics[width=2.75in]{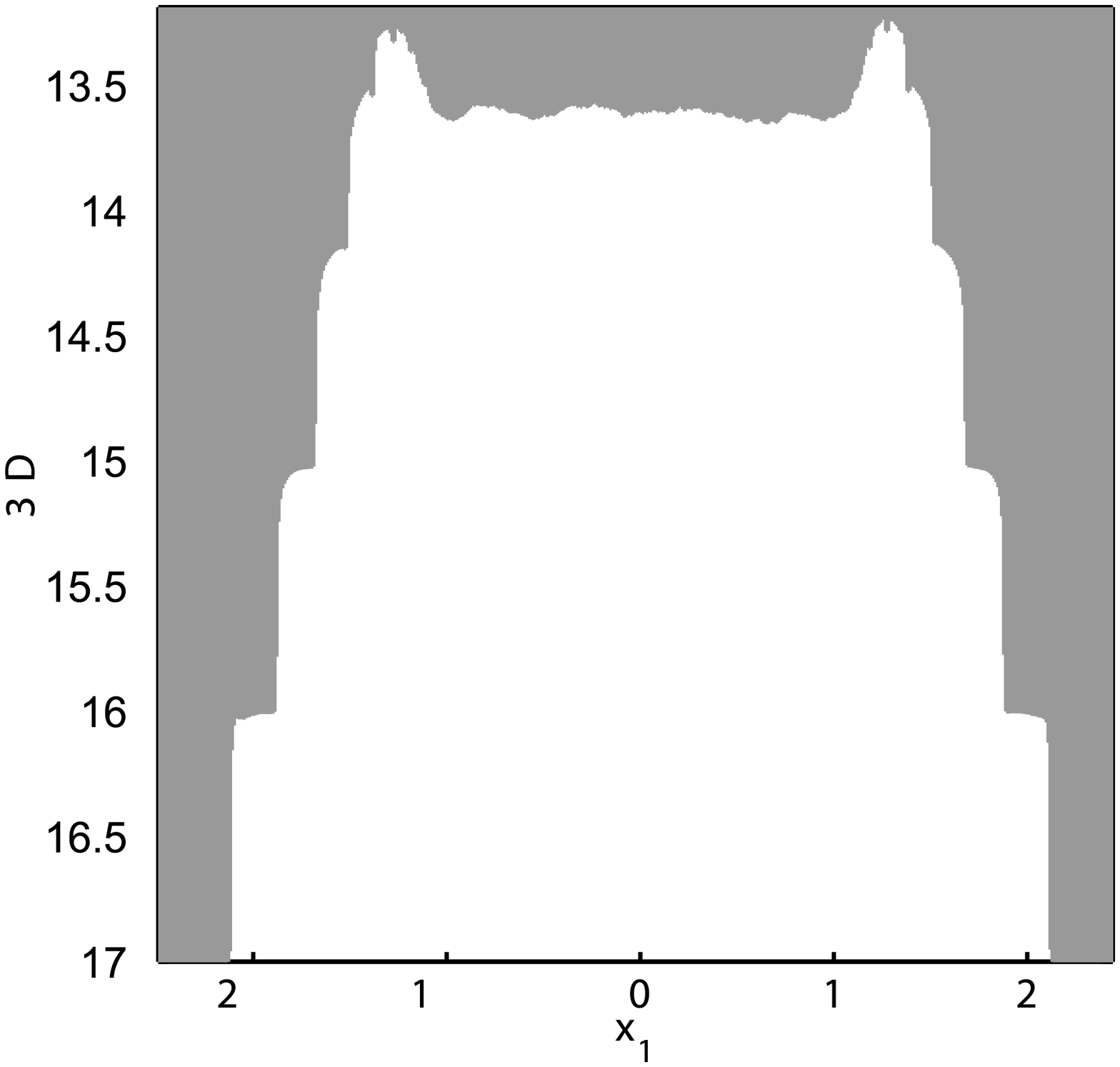} &  \includegraphics[width=2.75in]{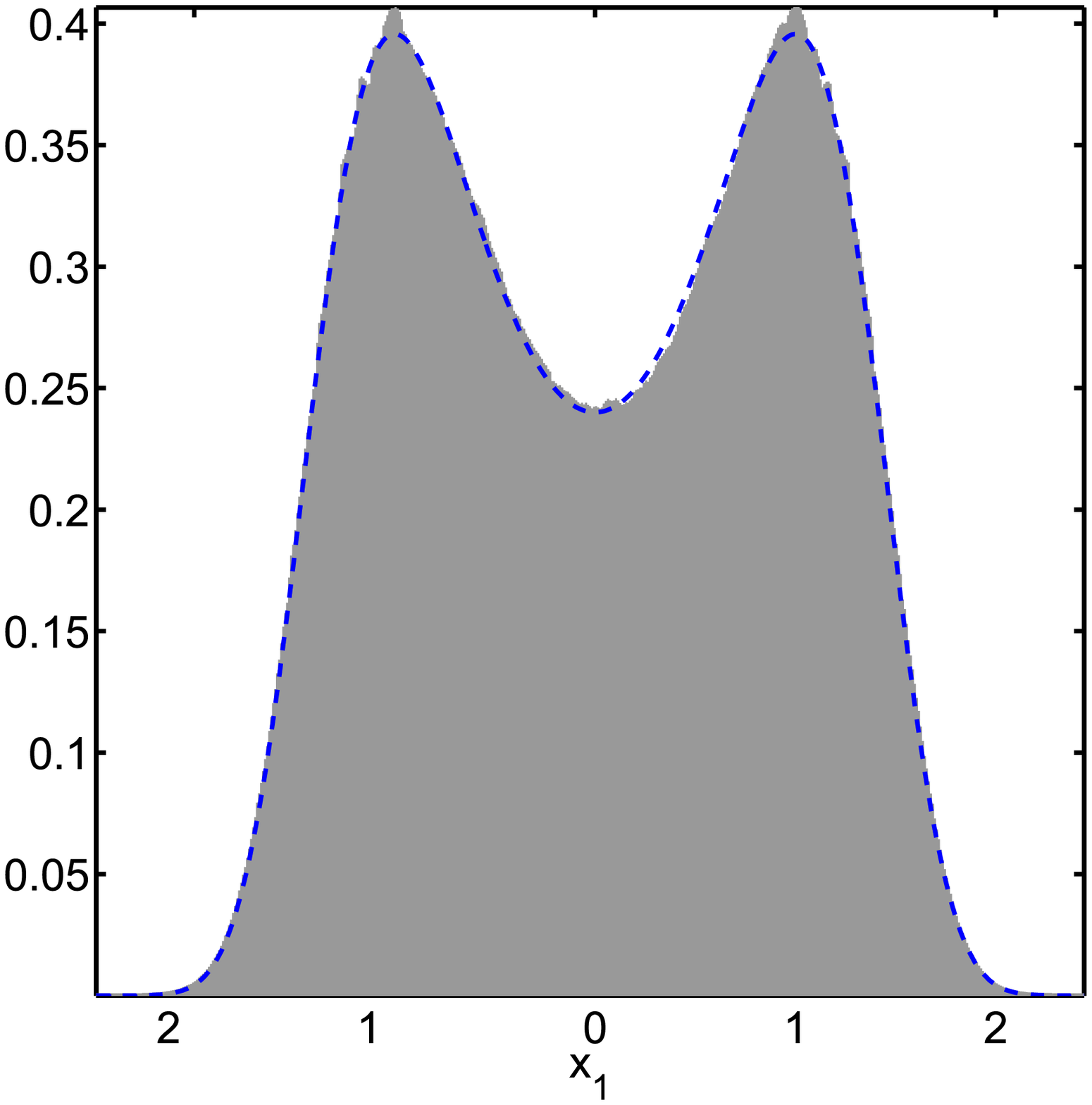}
  \end{tabular}
  \caption {Effect of the ``bad set'' $\Psi(\epsilon)$ in different
    dimensions.  On the right, the numerically obtained density and true
    analytic density (dashed line). On the left, The average $log$ of the step
    size taken verses the spatial position of $x_1$.}
 \label{fig:grad-problem}%
\end{figure}

\begin{figure}[bthp]
  \begin{center}
    \scalebox{0.7}{\includegraphics{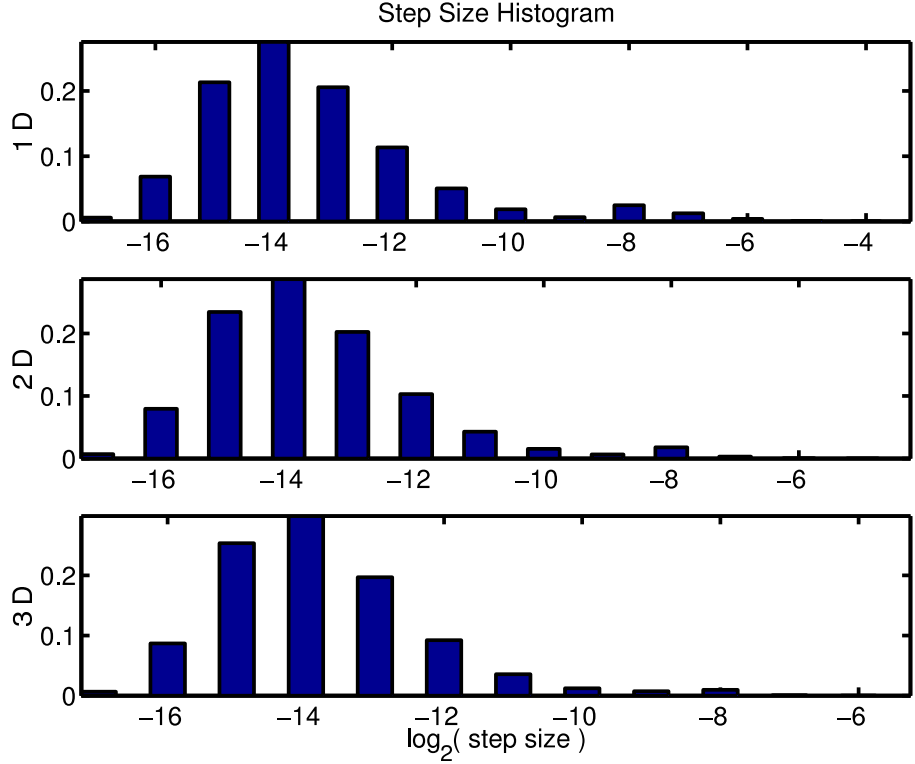}}
  \end{center}
  \caption {Gradient Problem -- Timestep Historgram}
  \label{fig:grad-steps}%
\end{figure}

For contrast we apply the algorithm to the systems in two and three
dimensions found by making identical copies of the equation
(\ref{eq:test2}) in the extra dimensions with each dimension driven by
an independent noise. Thus the set $\Psi(0)$ comprises the tensor
product of the set $\pm 1,0, \pm \frac{1}{\sqrt 3}$ in the appropriate
number of dimensions. Small neighbourhoods of this set do have large
hitting time, with high probability. To illustrate the effect of this
difference between hitting times we show, in Figure
\ref{fig:grad-problem}, the average time step taken at a given
location in space for (the first component of) $y$. Notice that in one
dimension the algorithm allows quite large average steps in the
neighbourhood of the points $\pm 1,0, \pm \frac{1}{\sqrt 3}.$ This
does not happen in dimensions two and three because the probability
that the other components of $y$ is also near to the set $\pm 1,0, \pm
\frac{1}{\sqrt 3}$ at the same time is very small. The effect of this
large choice of time steps in one dimension is apparent in the
empirical densities for (the first component of) $y$ which are also
shown in Figure \ref{fig:grad-problem}; these are generated by binning
two hundred paths of the SDE (\ref{eq:test2}) over two hundred time
units.  It is important to realize that, although the algorithm in one
dimension makes a very poor approximation of the empirical density,
this occurs only because of a relatively small number of poorly chosen
time-steps. Figure \ref{fig:grad-steps} shows a histogram of the
timesteps ($k_n$ values) taken in one, two and three dimensions. The
plots are nearly identical, except that in one dimension the algorithm
allows the method to take a small number of larger steps with $k_n=4.$

\section{Generalizations of the Method}
\label{sec:generalization}

The method given in \eqref{eq:ratz} can be seen as a simple
instance of a general class of methods based on comparing, with some error metric, one time step given by two methods
of the form
\begin{align}
  x_{n+1} &=  x_n + F(x_n,\Dt_{n+1}) \Dt_{n+1} + G(x_n,\Dt_{n+1})\sqrt \Dt_{n+1}\eta_n \label{eq:genForm}\\
  \overline{x}_{n+1} &= x_n +
  \overline{F}(x_n,\Dt_{n+1}) \Dt_{n+1} +
  \overline{G}(x_n,\Dt_{n+1})\sqrt
  \Dt_{n+1}\eta_n \;, \notag
\end{align}
where $F,\overline{F},G,\overline{G}$ are deterministic functions.
The method in \eqref{eq:ratz} was based on comparing the pair of
explicit methods given by
\begin{align}
  \label{eq:ratz-again}
  x_{n+1} &= \xs_n+\sqrt \Dt_{n+1}\Sigma(x_n)\eta_n\\
  \overline{x}_{n+1} &= \hat{x}_n+\sqrt \Dt_{n+1}\Sigma(\xs_n)\eta_n \notag
\end{align}
where
\begin{align*}
  \xs_n = & x_n +\Dt_{n+1} f(x_n),\\
  \hat{x}_n = & x_n+\frac12 \Dt_{n+1}[f(x_n)+f(\xs_n)] .  \end{align*}
In \eqref{eq:ratz}, closeness was measured by the difference, divided
by the step size, between the conditional expectations of one time step 
of the two different methods; this gives
\begin{align}\label{eq:drift_control}
\frac{2}{\Dt_{n+1}} |\EE x_{n+1} - \EE\overline{x}_{n+1}| =
2|F(x_n,\Dt_{n+1}) - \bar F(x_n,\Dt_{n+1})|=|f(x_n) - f(\xs_n)| \ .
\end{align}

From this point of view, it is clear that the method discussed thus
far is one of a large family of methods. Depending on the setting, one
might want to compare methods other that the simple Euler methods
used thus far. Also one can consider different error measures. In the
next section, we study a damped-driven Hamiltonian problem and use ideas
from symplectic integration to design an appropriate method. In the
discussion at the end of the article, we return to the question of
different error measures.

\section{Numerical Experiments: Long Time Simulations}
\label{sec:9}

In this section, we demonstrate that the ideas established for the rather
specific adaptive scheme studied, and for the particular hypotheses
on the drift and diffusion, extend to a wider class of SDEs and
adaptive methods. 

As a test problem we consider the Langevin equation
\begin{align}\label{eq:Langevin}
  dq=&p\; dt\\
  dp=& -\big[\delta(q)p+\Phi'(q)\big]\;dt+\Sigma(q)dW \ . \notag
\end{align} 
where $2\delta(q)=\Sigma^2(q)$,
\begin{equation*}
  \Phi(q)=\frac14(1-q^2)^2\quad\mbox{and}\quad 
  \Sigma(q)=\frac{4(5q^2 + 1)}{5(q^2 +1)} .
\end{equation*}
The preceding theory does not apply to this system since it
is not uniformly elliptic; furthermore it fails to satisfy \eqref{eq:coerce}. 
However it does satisfy a Foster-Lyapunov drift condition and since it
is hypoelliptic the equation itself can be proven geometrically
ergodic \cite{msh}. In \cite{msh}, it was shown that the implicit
Euler scheme was ergodic when applied to \eqref{eq:Langevin}, and a similar
analysis would apply to a variety of implicit methods. Since
the adaptive schemes we study in this section
enforce closeness to such implicit methods, we believe that analysis 
similar to that in the previous section will extend to this Langevin equation
and to the adaptive numerical methods studied here.

We will compare two different methods based on different choices  of
the stepping method. The first is the Euler
based scheme given in \eqref{eq:ratz}. The second is the following 
scheme:
\begin{align}
  \label{symplectic}
  q_{n+1} = & \qs_n \\
  p_{n+1} = & \ps_n + \Sigma(\qs_n)\sqrt{\Dt_n}\eta_{n+1}  \notag \\
  \overline{q}_{n+1} =& q_n +\Big(\frac{p_n + \ps_n}{2}\Big) \Dt_n
  \notag 
  \\
  \overline{p}_{n+1} = & p_n -\Phi'\Big(\frac{q_n + \qs_n}{2}\Big) \Dt_n -
  \delta\Big(\frac{q_n + \qs_n}{2}\Big)\Big(\frac{p_n
  +\ps_n}{2}\Big)\Dt_n + \Sigma\Big(\frac{q_n
  +\qs_n}{2}\Big)\sqrt{\Dt_n}\eta_{n+1}  \notag
\end{align}
where
\begin{align*}
  \qs_n=& q_n + p_n\Dt_n,\\
  \ps_n=& p_n - \Phi'(\qs_n)\Dt_n - \delta(\qs_n)p_n \Dt_n.
\end{align*} Once again we will use comparisons between the two updates
(with and without bars) to control the error. As before, we control on
the difference in the expected step. The point of the particular form
used here is that, in the absence of noise and damping, the adaptation
constrains the scheme to take steps which are close to those of the
symplectic midpoint scheme, known to be advantageous for Hamiltonian
problems; if the noise and damping are small, we expect the
Hamiltonian structure to be important.

In both the Euler and Symplectic case,  the stepping methods take the form 
\begin{align}\label{eq:methods}
  q_{n+1}&=q_n+f_n^{(1)}\Dt_n\\
  p_{n+1}&=p_n+f_n^{(2)}\Dt_n + g_n\sqrt{\Dt_n}\eta_{n+1} \notag\\
  \overline{q}_{n+1}&=\overline{q}_n+\overline{f}_n^{(1)}\Dt_n \notag\\
  \overline{p}_{n+1}&=\overline{p}_n+\overline{f}_n^{(2)}\Dt_n +
  \overline{g}_n\sqrt{\Dt_n}\eta_{n+1} \notag
\end{align}
where $f_n$ and $\overline{f}_n$, $g_n$, $\overline{g}_n$ are
adapted to $\cF_{n-1}$. In this notation the metric becomes,
\begin{align*}
 \big[ (\overline{f}_n^{(1)} - f_n^{(1)})^2 +  (\overline{f}_n^{(2)} -
 f_n^{(2)})^2 \big]^\frac12 < \tau
\end{align*}

In the remainder of this section, we present numerical experiments
with the two methods just outlined. We study the qualitative approximation of
the invariant measure, we quantify this approximation and measure its
efficiency, and we study the behaviour of time-steps generated.

 \begin{figure}[tbhp]\label{fig:invMeasureEuler1}
   \centering
   \includegraphics[width=3in]{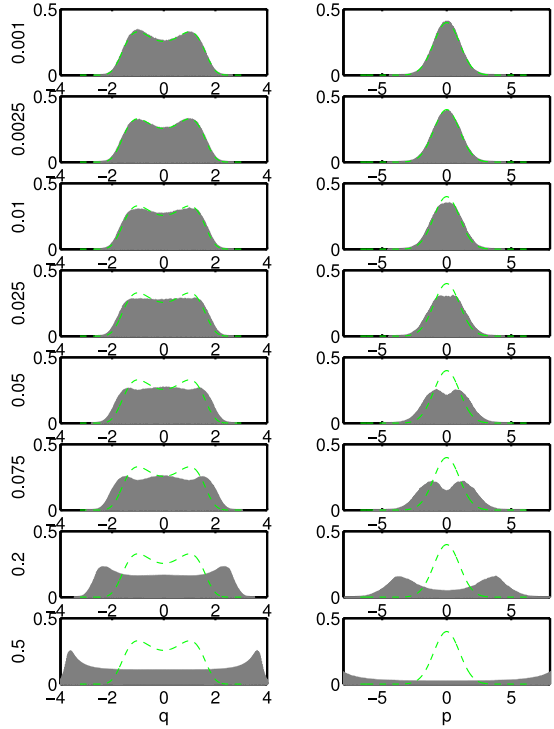} 
   \caption{Distribution of $q,p$ with different $\tau$ for Euler 
       method The value of tolerance
       $\tau$ is on the left of each figure.}
   \protect\label{fig:pdfEuler1}
 \end{figure}

 \begin{figure}[tbhp]\label{fig:invMeasureEuler2}
   \centering
       \includegraphics[width=3in]{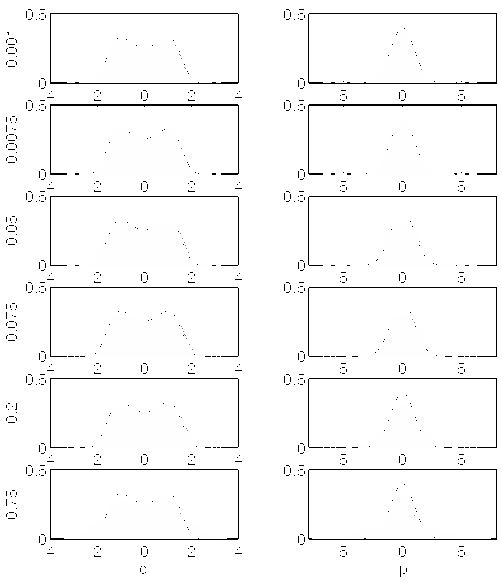} 
   \caption{Distribution of $q,p$ with different $\tau$ for the symplectic 
       method. The value of tolerance
       $\tau$ is on the left of each figure.}
   \protect\label{fig:pdfEuler2}
 \end{figure}

 Figure \ref{fig:invMeasureEuler1} plots the numerically obtained time
 average of the position $p$ and momentum $q$ for various values of
 the tolerance $\tau$. The doted lines are the true invariant measure
 of the underling SDE which can be computed analytically in this
 particular case. Notice that the method appears stable for all values
 of $\tau$, in contrast to the forward Euler method which blows up
 when applied to this equation. Though apparently stable, the results
 are far from the true distribution for large $\tau$.
 Figure \ref{fig:invMeasureEuler2} gives the analogous plots for the
 adaptive symplectic method given in \eqref{eq:methods}.  Notice that
 these methods seem to do a much better job of reproducing the invariant
measure faithfully at large $\tau$.

\begin{figure}[htbp]
   \centering
    \includegraphics[height=5in,angle=-90]{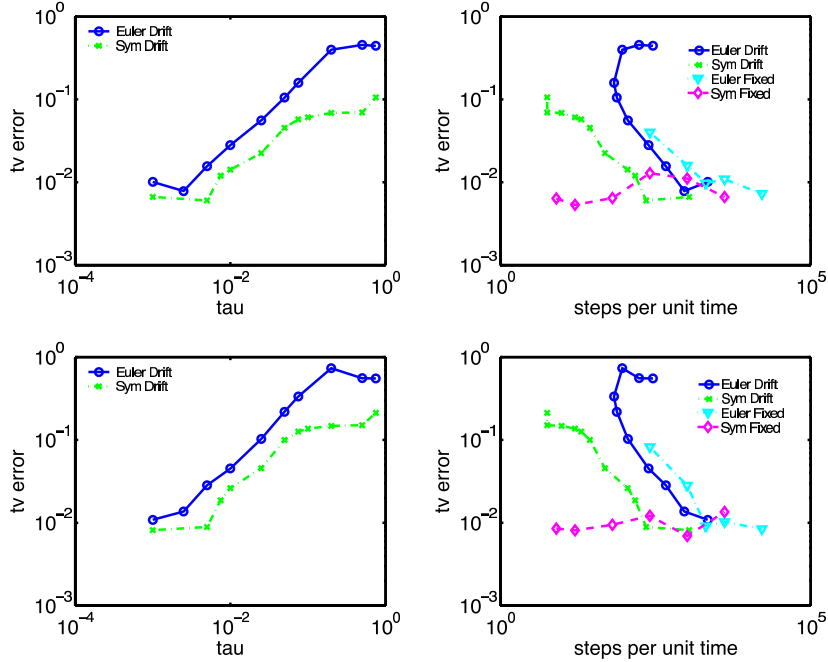}
   \caption{ Total variation error verses $\tau$ for position (top
    left) and the momentum (bottom left). Total variation error verses
   $\tau$ for steps per unit time (top right) and the momentum (bottom right). }
   \protect\label{fig:compare}
 \end{figure}

 It is also important to study accuracy per unit  
 of computational effort. Figure \ref{fig:compare} gives plots
 of the error in the total variation norm (the $L^1$ distance between
 the numerically computed time averages and the exact analytic answer
 verses the $\tau$ used and versus the steps per unit of time; the latter
provides a measure of unit cost. The top
 plots are for the momentum $q$ and the bottom for the position $p$.
 The plots on the right also include two fixed step methods, one using
 simple the forward Euler scheme and the second using the first of the
 symplectic schemes. The fixed-step Euler schemes blows up for steps
 larger than those given. We make the following observations on the
basis of this experiment:

\begin{itemize}

\item The fixed-step symplectic method is the most efficient at small time-steps;

\item The  adaptive symplectic method is considerably more efficient than
the adaptive and fixed-step Euler methods;

\item The adaptive symplectic method is the most robust method, providing
reasonable approximations to the invariant density for a wide range of $\tau.$

\end{itemize}

Note that the adaptive methods have not been optimized and with careful
attention might well beat the fixed-step methods, both as measured
by accuracy per unit cost, as well as by robustness. Further 
study of this issue is required.

 \section{Conclusions}
 \label{sec:conc}
 \setcounter{equation}{0}

This paper proposes a simple adaptive strategy for SDEs which is
designed to enforce geometric ergodicity, when it is present
in the equation to be approximated; without adaptation
methods such as explicit Euler may destroy ergodicity. As well
as proving ergodicity, we also prove some exponential moment
bounds on the numerical solution, again mimicking those for the SDE
itself. Furthermore, we prove finite-time convergence of the numerical method;
this is non-trivial because we do not assume (and it is not true in
general) that the time-step sequence tends to zero with user input tolerance.
It would be of interest to transfer this finite time convergence to a result
concerning convergence of the invariant measures, something which is known
for fixed time-step schemes \cite{Tal,Talnew,Talnew2}.

As discussed in section \ref{sec:generalization}, the scheme we study
in detail here is prototypical of more advanced schemes comparing two
more sophisticated methods and controlling both on drift and
diffusion.  Here we have mainly used simple forward Euler methods and
controlled only on the drift: our error measure is based on the conditional
means. 
The split-step approach we take allows for additional
terms to be added to the error measure, to ensure that the diffusion
step is also controlled.  The general idea is to enforce the closeness
of one step by two different methods. One has freedom in the choice of
the methods and the measure of closeness. We now briefly mention to
other error measure which make this idea specific.

For simplicity let us assume work in one dimension though the ideas
generalize directly to higher dimensions. The simple error control
given in \eqref{eq:drift_control} controls only the difference in the
expectation of one step of the two methods. However one can also use
measure which ensure the closeness of the entire distribution of one time step
of the two methods. Given $x_n=\bar x_n$ and $\Delta_{n+1}$, one
step of a method of the form  \eqref{eq:genForm} is Gaussian. Hence it is
reasonable to require that the standard deviations are close to
each other. The error criteria would then be 
\begin{multline*}
  \frac{1}{\Dt_{n+1}} |\EE   x_{n+1} - \EE \overline{x}_{n+1}|+  {\frac{1}{ \sqrt{\Dt_{n+1}}} |
\mbox{StdDev}(x_{n+1})  - \mbox{StdDev}(\bar x_{n+1})|} =\\|F(x_n,\Dt_{n+1}) -
\bar F(x_n,\Dt_{n+1})| + \big| |G(x_n,\Dt_{n+1})| - |\bar
G(x_n,\Dt_{n+1})| \big| < \tau
\end{multline*}
In some ways, comparing the mean and standard deviations is rather arbitrary.
A more rational choice might be to ensure the closeness of the total
variation distance  of the densities after one time step of the two
methods. A simple way to do this is to compare the relative entropy of
the two distributions. Since the distributions are Gaussian this can
be done explicitly. One finds that the criteria based on controlling relative entropy per unit step
is
\begin{align*}
  \frac{\big(F(x_n,\Dt_{n+1}) -
\bar F(x_n,\Dt_{n+1}))^2}{G(x_n,\Dt_{n+1}\big)^2} +\Big(\frac{\bar
G(x_n,\Dt_{n+1})^2}{G(x_n,\Dt_{n+1})^2} -1\Big) -\log \frac{\bar
G(x_n,\Dt_{n+1})^2}{G(x_n,\Dt_{n+1})^2} < 2 \tau \ .
\end{align*}
It is interesting to note that this measure correctly captures the
fact that one should measure the error in the drift on the scale of
the variance. In other words if the variance is large, one does not
need to be as accurate in calculating the drift as it will be washed
out by the noise anyway. Since the above measure is expensive to
calculate one can use that fact that $\frac{\bar G}{G} -1$ is small to
obtain the asymptotically equivalent criteria
\begin{align*}
  \frac{\big(F(x_n,\Dt_{n+1}) -
\bar F(x_n,\Dt_{n+1}))^2}{G(x_n,\Dt_{n+1}\big)^2} +\frac12\Big(\frac{\bar
G(x_n,\Dt_{n+1})^2}{G(x_n,\Dt_{n+1})^2} -1\Big)^2 < 2 \tau \ .
\end{align*}

\section{Acknowledgments}

The authors thank George Papanicolaou for useful discussions
concerning this work. JCM thanks the NSF for its support (grants
DMS-9971087 and DMS-9729992); he also thanks the Institute for
Advanced Study, Princeton, for its support and hospitality during the
academic year 2002-2003. AMS thanks the EPSRC for financial support.
We also wish to thank a careful and conscientious referee who noticed a
substantive error in an earlier draft of this paper, and persisted until we
understood his/her point.


\appendix 

\section{Two Exponential Martingale Estimates in Discrete Time}

Let $\{\mathcal{F}_n$, $n \geq 0\}$ be a filtration. Let $\eta_k$ be a
sequence of random variables with $\eta_k$ adapted to $\mathcal{F}_k$
and such that $\eta_{k+1}$ conditioned on $\mathcal{F}_{k}$ is normal
with mean zero and variance $\sigma_{k}^2=\EE[
\eta_{k+1}^2|\mathcal{F}_k]< \infty$. We define the following processes:
\begin{xalignat*}{2}
  M_n&=\sum_{k=1}^n \eta_k &  \tilde M_n&=\sum_{k=1}^n \eta_k^2 -
  \sigma_{k-1}^2\\
  \la M\ra_n&= \sum_{k=0}^{n-1} \sigma_k^2 &   \la \tilde M\ra_n&= \sum_{k=0}^{n-1} 2\sigma_k^4 
\end{xalignat*}
 As the notation suggests, $ \la M\ra_n$ and $ \la \tilde
M\ra_n$ are the quadratic variation processes in that $M_n^2-\la
M\ra_n$ and $\tilde M_n^2-\la \tilde M\ra_n$ are local martingales with
respect to $\mathcal{F}_n$.

\begin{lemma}
  \label{l:expMartNew} Let $\alpha >0$ and $\beta>0$, then the
  following estimate holds:
  \begin{align*}
    \PP\Big( \sup_{k} \big(  M_k -  \frac\alpha2\la M
    \ra_k\big) \geq  \beta \Big) & \leq e^{-\alpha\beta}\\
\intertext{If in addition $\sigma_k^2 \leq \sigma_*^2 \in \RR$ for all $k
\in\NN$ almost surely, then}
    \PP\Big( \sup_{k} \big( \tilde M_k -\frac\alpha2 \la \tilde M
    \ra_k\big) \geq \beta \Big) &\leq e^{-\frac{\beta}{\lambda^2}}
  \end{align*}
where $\lambda^2= 2 \sigma_*^2 + 1/\alpha$.
\end{lemma}

\begin{proof}[Proof of Lemma \ref{l:expMartNew}]
  We begin with the first estimate. Define $N_n=\exp( \alpha M_n -
  \frac{\alpha^2}2\la M \ra_n)$ and observe that $N_{n} = \EE\{N_{n+1}
  | \mathcal{F}_n\}$.  This in turn implies that $\EE|N_n|=\EE N_n =
  N_0=1 < \infty$. Combining these facts we see that $N_n$ is a
  Martingale. Hence, the Doob-Kolmogorov Martingale inequality
  \cite{RW} implies
  \begin{align*}
    \PP\big( \sup_n N_n > c \big) \leq \frac{\EE N_0}{c} = \frac1c \ .
  \end{align*}
  Since $\PP\Bigl( \sup_n (M_n - \frac{\alpha}2\la M\ra_n) \ge  \beta \Bigr)=
  \PP\big( \sup_n N_n > e^{\alpha\beta} )$, the proof is complete.

  The second estimate is obtained in the same way after some
  preliminary calculations. We define $\phi(x)= \frac12 \ln( 1- 2 x)$
  and $\psi(x,b)=-x -b x^2$. Observe that $ c \psi(x,b)=\psi(c x,b/c)$
  and $\phi(x) \geq \psi(x,b)$ if $ x \in \big[ 0,
  \frac12\big(\frac{b-1}{b}\big)\big]$ and $b >1$. Now
  \begin{align*}
    \PP\left( \sup_n (\tilde M_n - \frac{\alpha}{2}\langle\tilde M
      \rangle_n)  \ge  \beta \right) & = \PP\left( \sup_n \sum_{k=1}^n
      \frac{\eta_k^2}{\lambda^2} +
      \frac{1}{\lambda^2}\psi(\sigma_{k-1}^2,\alpha) \ge
      \frac\beta{\lambda^2} \right)\ .
  \end{align*}
  Setting $\lambda^2= 2 \sigma_*^2 + \frac1\alpha$, we have that
  $\frac{1}{\lambda^2}\psi(\sigma_{k}^2,\alpha)=
  \psi(\frac{\sigma_{k}^2}{\lambda^2},\lambda^2\alpha) \leq
  \phi(\frac{\sigma_{k}^2}{\lambda^2})$ for all $k \geq 0$ since
  $\sigma_k^2 \leq \sigma_*^2$ and $\lambda^2 \alpha > 1 $. Defining 
  \begin{align*}
    \tilde N_n = \exp\left( \sum_{k=1}^n \frac{\eta_k^2}{\lambda^2} +
    \phi(\frac{\sigma_{k-1}^2}{\lambda^2})\right)\ ,
  \end{align*}
  we have
  \begin{align*}
    \PP\left( \sup_n (\tilde M_n - \frac{\alpha}{2}\langle\tilde M
      \rangle_n) \ge \beta \right) &\leq \PP\left( \sup_n \tilde N_n \ge
      e^\frac{\beta}{\lambda^2} \right)\ .
  \end{align*}
  Now recall that if $\xi$ is a unit Gaussian random variable then
  $\EE \exp( c \xi^2)= 1/\sqrt{1-2c}$ for $c \in (-\frac12,\frac12)$.
  By construction $\frac{\eta_k}{\lambda}$, conditioned on
  $\mathcal{F}_{k-1}$, is a Gaussian random variable with variance
  less then $\frac12$. Hence
  \begin{equation*}
\EE\Bigl( \exp\big(\frac{\eta_k^2}{\lambda^2}\big)\big| \mathcal{F}_{k-1}
\Bigr) = \exp\Bigl(-\phi\bigl(\frac{\sigma_{k-1}^2}{\lambda^2}\bigr)\Bigr)\;.
  \end{equation*}
  Using this one sees that $\EE \big\{\tilde N_{n+1} |
  \mathcal{F}_n\big\} = \tilde N_n$ and $\EE |\tilde N_n|=1 < \infty$,
  hence $\tilde N_n$ is a Martingale.   By the same argument
  as before using the Doob-Kolmogorov Martingale inequality, we obtain
  the quoted result.
\end{proof}


{ \small
\begin{thebibliography}{99}








\bibitem{AGH} M.A. Aves, D.F. Griffiths and D.J. Higham,
{\it Does error control suppress spuriosity?}
To appear in SIAM J. Num. Anal.

\bibitem{CR} G. Cassella and C.P. Robert, {\it Monte Carlo
Statistical Methods}, Springer Texts in Statistics, New York,
2002.

\bibitem{DV} K. Dekker and J.G Verwer, {\it Stability of
Runge-Kutta Methods for Stiff Nonlinear Differential Equations.}
North-Holland, Amsterdam, 1984.

\bibitem{Dud} R.M. Dudley, {\it Real Analysis and Probability}.
Cambridge University Press, Cambridge, 2002.

\bibitem{GL} J.G. Gaines and T.J. Lyons,
{\em Variable stepsize control in the numerical solution of stochastic
differential equations}, SIAM J. Appl. Math., {\bf 57}(1997),1455--1484.

\bibitem{Gr} D.F. Griffiths, {\it The dynamics of some linear
multistep methods with step-size control.} Appears in
``Numerical Analysis 1987''
Eds: Griffiths, D.F. and Watson, G.A., Longman (1988), 115--134.

\bibitem{Ha} J.K. Hale, {\it Asymptotic Behaviour
of Dissipative Systems.} AMS Mathematical Surveys and
Monographs 25, Rhode Island, 1988.

\bibitem{hansen} N.R. Hansen, {\em Geometric ergodicity of
discrete time approximations to multivariante diffusions}.
Bernoulli, submitted, 2002.


\bibitem{has} 
R.~Z. {Has'minski{i}}.
\newblock {\em Stochastic Stability of Differential Equations}.
\newblock Sijthoff and Noordhoff, 1980.


\bibitem{h&s}  D.J. Higham and A.M. Stuart
{\em Analysis of the dynamics of error control via a piecewise
continuous residual.} BIT {\bf 38}(1998), 44--57.


\bibitem{HMS}
D.~J. Higham, X.~Mao, and A.~M. Stuart.
\newblock Strong convergence of {E}uler-type methods for nonlinear stochastic
  differential equations.
\newblock {\em SIAM J. Numer. Anal.}, 40:1041--1063, 2002.

\bibitem{ritter} N. Hofmann, T. M\"{u}ller-Gronbach and K. Ritter
{\em Optimal approximation of SDEs by adaptive step-size control.}
Math. Comp. {\bf 69}(2000), 1017--1034.


\bibitem{ritter2} N. Hofmann, T. M\"{u}ller-Gronbach and K. Ritter
{\em The optimal discretization of SDEs.}
J. of Complexity. {\bf 17}(2001), 117--153.

\bibitem{ritter3} N. Hofmann, T. M\"{u}ller-Gronbach and K. Ritter
{\em Optimal uniform approximation of systems of SDEs.}
Ann. Appl. Prob. {\bf 12}(2002), 664--690.

\bibitem{KP} P.E. Kloeden and E. Platen, {\em Numerical Solution of
Stochastic Differential Equations}. Springer-Verlag, New York, 1991.


\bibitem{L} H.~Lamba.
\newblock Dynamical systems and adaptive time-stepping in {ODE} solvers.
\newblock {\em BIT}, 40:314--335, 2000.



\bibitem{HL} H. Lamba and A.M. Stuart {\em Convergence results
for the MATLAB ode23 Routine}. BIT {\bf 38}(1998), 751--780.

\bibitem{LMS} H. Lamba, J.C. Mattingly and A.M. Stuart {\em 
An adaptive Euler-Maruyama scheme for SDEs: Part I, Convergence}.

\bibitem{mao} X. Mao, {\em Stochastic Differential Equations and Applications}.
Horwood, Chichester, 1997.

\bibitem{Matlab}
  The MathWorks, Inc.,
  {\em MATLAB User's Guide}.
  Natick, Massachusetts,
  1992.

\bibitem{msh} J.C.  Mattingly, A.M. Stuart and D.J. Higham
{\em Ergodicity for SDEs and approximations: locally Lipschitz
vector fields and degenerate noise}.  Stoch. Proc. and Applics.
{\bf 101}(2002), 185--232.

\bibitem{mat} J.C. Mattingly, {\em A numerical study of adaptive
methods for SDEs}. In preparation.

\bibitem{MT} S. Meyn and R.L. Tweedie, {\em Stochastic Stability
of Markov Chains}. Springer-Verlag, New York, 1992.

\bibitem{MG2}
T.~M\"{u}ller-Gronbach.
\newblock {\em Strong approximation of systems of stochastic differential
  equations}.


\bibitem{RT} G.O. Roberts and R.L. Tweedie
{\em Exponential convergence of Langevin diffusions and their discrete
approximations},  Bernoulli {\bf 2}(1996), 341--363.








\bibitem{RW} L.C.G. Rogers and D.Williams {\em Diffusions, Markov
processes and Martingales, Volumes 1 and 2}. Cambridge University Press,
reprinted second edition, 2000.


\bibitem{Sanz} J.-M. Sanz-Serna, {\em Numerical ordinary differential
equations versus dynamical systems}. In D.S. Broomhead and A. Iserles,
Editors, ''The Dynamics of Numerics and the Numerics of Dynamics'',
Clarendon Press, Oxford, 1992.


\bibitem{s&t} O. Stramer and R.L. Tweedie, {\em Langevin-type models I:
diffusions with given stationary distributions, and their discretizations},
Methodology \& Computing in Applied Probability, {\bf 1}(1999), 283-306. 

\bibitem{S}
A.M. Stuart.
\newblock Probabilistic and deterministic convergence proofs for software for
  initial value problems.
\newblock {\em Numerical Algorithms}, 14:227--260, 1997.

\bibitem{SH} A.M. Stuart and A.R. Humphries,  {\it
Dynamical Systems and Numerical Analysis}.
Cambridge University Press, 1996.

\bibitem{SH2} A.M. Stuart and A.R. Humphries,  {\it
The essential stability of local error control for
dynamical systems.}
SIAM J. Num. Anal. {\bf 32}(1995), 1940--1971.


\bibitem{Tal} D. Talay, {\em Second-order discretization schemes
for stochastic differential systems for the computation of the
invariant law.} Stochastics and Stochastics Reports {\bf 29}(1990),
13--36.

\bibitem{Talnew} D. Talay, {\em Approximation of the invariant probability
measure of stochastic Hamiltonian dissipative systems with non globally
Lipschitz co-efficients}. Appears in ``Progress in Stochastic
Structural Dynamics'', Volume 152, 1999, Editors R. Bouc
and C. Soize. Publication du L.M.A.-CNRS.

\bibitem{Talnew2} D. Talay, {\em
Stochastic Hamiltonian dissipative systems with non globally
Lipschitz co-efficients: exponential convergence to the invariant
measure and discretization by the implicit Euler scheme}.
Markov Proc. Rel. Fields., {\bf 8}(2002), 163--198. 

\bibitem{Te} R. Temam, {\it Infinite Dimensional Dynamical
Systems in Mechanics and Physics. Springer, New York, 1989.}






\end{thebibliography}
}
\end{document}